\documentclass[11pt]{amsart}
\usepackage[a4paper,top=4cm,bottom=3cm,left=3cm,right=3cm]{geometry}
\usepackage[T1]{fontenc}
\usepackage{bm}
\usepackage{amsmath,amssymb,amsfonts,amsthm} 
\usepackage{hyperref}
\usepackage{graphicx}
\usepackage{dirtytalk}
\usepackage{pgfplots}
\pgfplotsset{compat=newest}
\usepackage{yfonts}
\usepackage{url}
\usepackage{enumitem}
\usepackage{esint}
\usepackage{multirow}
\usepackage{blindtext}
\newcommand{\numberset}{\mathbb}
\newcommand{\N}{\numberset{N}}
\newcommand{\R}{\numberset{R}}
\newcommand{\M}{\textup M}
\newcommand{\g}{\textup g}
\newcommand{\K}{\textup K}
\newcommand{\J}{\textup J}

\newcommand{\cs}{\textup{c}_\textup{s}}
\newcommand{\cper}{\textup{c}_\textup{Per}}

\newcommand{\cin}{\textup{c}_{\textup{in}}}

\newcommand{\Y}{\Upsilon}
\newtheorem{thm}{Theorem}[section]
\newtheorem{lem}[thm]{Lemma}
\newtheorem{dfn}[thm]{Definition}
\newtheorem{cor}[thm]{Corollary}
\newtheorem{prp}[thm]{Proposition}
\newtheorem{rem}[thm]{Remark}

\title[Volume preserving spacetime mean curvature flow]{Volume preserving spacetime mean curvature flow and foliations of initial data sets}

\author{Jacopo Tenan}
\address{(Jacopo Tenan) Dipartimento di Matematica, Universit\`a degli Studi di Roma "Tor Vergata", Via della Ricerca Scientifica, 00133, Roma, Italy}
\date{}
\numberwithin{equation}{section}
\begin{document}
\sloppy
\begin{abstract}
We consider a volume preserving curvature evolution of surfaces in an asymptotically Euclidean initial data set with positive ADM-energy. The speed is given by a nonlinear function of the mean curvature which generalizes the spacetime mean curvature recently considered by Cederbaum-Sakovich (Calc. Var. PDE, 2021). Following a classical approach by Huisken-Yau (Invent. Math., 1996), we show that the flow starting from suitably round initial surfaces exists for all times and converges to a constant (spacetime) curvature limit. This provides an alternative construction of the CSTMC foliation by Cederbaum-Sakovich and has applications in the definition of center of mass of an isolated system in General Relativity.
\end{abstract}
\maketitle
\section{Introduction}
\noindent We consider a triple $(\M, \overline \g,\overline\K)$, where $(\M,\overline\g)$ is a smooth complete three-dimensional Riemannian manifold and $\overline\K$ a symmetric $(0,2)$-tensor field on $\M$, such that
\begin{equation}\label{einstconstsyst}
\begin{cases}
\overline{\textup{S}}-|\overline \K|^2+(\textup{tr}_{\overline \g}\overline \K)^2=2\overline\mu\\
\overline\nabla\cdot\left(\overline \K\right)-\overline {\textup d}\left(\textup{tr}_{\overline \g}\overline \K\right)=\overline \J
\end{cases}
\end{equation}
holds for a smooth function $\overline\mu\in C^\infty(\M)$ and a smooth 1-form $\overline\J \in \Omega^\infty(\M)$. For our purposes, we will assume that $(\M,\overline\g,\overline\K)$ and $(\overline\mu,\overline\J)$ satisfy the following definitions.
\begin{dfn}\label{def1} Let $\delta\in(0,\frac12]$. A $C_{\frac12+\delta}^2$-\textup{asymptotically flat initial data set} is a triple $(\M,\overline \g,\overline \K)$ such that there exist a compact subset $\textup C\subset \M$, a constant $\overline c>0$ and a diffeomorphism $\vec x:\M\setminus \textup C\to \R^3\setminus \overline{\mathbb B}_1(\vec 0)$ such that 
\begin{equation}\label{asymptc8}
|\overline\g_{\alpha\beta}-\delta_{\alpha\beta}|+|\vec x|\left|\partial_\gamma \overline\g_{\alpha\beta}\right|+|\vec x|^2\left|\partial_\gamma\partial_\epsilon\overline\g_{\alpha\beta}\right|\leq \overline c|\vec x|^{-\frac12-\delta},
\end{equation}
\begin{equation}\label{einstconstsyst2}
|\overline \K_{\alpha\beta}|+|\vec x||\partial_\gamma\overline \K_{\alpha\beta}|\leq \overline c|\vec x|^{-\frac32-\delta},
\end{equation}
where $\overline\g_{\alpha\beta}=\left(\vec x^*\overline \g\right)_{\alpha\beta}$ and $\overline \K_{\alpha\beta}=(\vec x^*\overline \K)_{\alpha\beta}$. Here $\partial_\gamma,\partial_\epsilon$ are partial derivatives with respect to the Euclidean coordinates, and the norms are taken in the Euclidean metric.
\end{dfn}
\begin{dfn}\label{def2}
We say that $(\M,\overline\g,\overline\K)$ is \textup{constrained} by the pair $(\overline\mu,\overline\J)$ if \eqref{einstconstsyst} holds together with 
\begin{equation}\label{constrainedmuJ}
|\overline\mu|+|\overline \J|\leq \overline c|\vec x|^{-3-\delta}.
\end{equation}
\end{dfn}
These triples $(\M,\overline\g,\overline\K)$ arise in Mathematical General Relativity as spacelike hypersurfaces of a Lorentzian spacetime $(\M,\overline \g)$ modelling an isolated gravitating system, with $\overline\K$ the second fundamental form of $\M$. In this description, it is important to give an appropriate definition of the center of mass of the system which is physically meaningful. A classical definition, inspired by the one of ADM-mass, was provided by Beig-Ó Murchadha \cite{beigomurchadha} as the limit of a flux integral on Euclidean spheres with increasing radius. A different approach was introduced by Huisken and Yau \cite{huiskenyau}, who proved the existence of a unique constant mean curvature (CMC)-foliation of the outer part of M, and defined the center of mass as the limit of the Euclidean barycenter of the leaves of the foliation as the radius tends to infinity. The result in \cite{huiskenyau} considered manifolds $\M$ which are asymptotically Schwarzschild and with strong decay assumptions. The construction was then generalized by various authors (e.g. \cite{ye}, \cite{metzger}, \cite{huang}) culminating in the work of Nerz \cite{nerz1} who considered asymptotically flat spaces with optimal decay assumptions. Under suitable symmetry assumptions, known as \textit{strong Regge-Teitelboim conditions}, it was proved that the CMC and the Beig-Ó Murchadha centers of mass both exist and coincide. On the other hand, on an asymptotically flat manifold where the (weak) Regge-Teitelboim condition are not satisfied,  these centers of mass may not be well-defined. In particular,  Cederbaum and Nerz \cite{cederbaumnerz} constructed explicit examples where both these objects do not converge.

For this reason, Cederbaum and Sakovich \cite{cederbaum} have introduced a new notion of center of mass based on foliations involving a modified curvature. For a smooth surface $\Sigma\hookrightarrow \M$, let us denote by $H$ the usual mean curvature, by $g=\iota^*\overline\g$ the induced metric on $\Sigma$, and by $P=\textup{tr}_g(\overline\K)$ the trace of $\overline\K$ on $\Sigma$. We then call the quantity $\mathcal H=\sqrt{H^2-P^2}$ the \textit{spacetime mean curvature} of $\Sigma$. We observe that $\mathcal H$ equals the (Minkowski) length of the mean curvature vector of the 2-codimensional surface $\Sigma$. In \cite{cederbaum}, Cederbaum and Sakovich prove the existence of a \textit{constant spacetime mean curvature} (CSTMC) \textit{foliation} of the outer part of $\M$, and define a corresponding CSTMC-center of mass as the limit of the barycenters of the leaves. It is proved that the CSTMC-center of mass exists in some relevant cases in which the CMC-center of mass does not, see \cite[Sect. 9]{cederbaum}.\\
\indent From a spacelike point of view, the equation satisfied by each CSTMC surface looks like a \textit{prescribed mean curvature equation}. A similar equation was present in the work of Metzger \cite{metzger}, who constructed surfaces satisfying the so called \textit{constant expansion equations} $\Theta_\pm := H \pm P\equiv\textup{const}$, where $\Theta_\pm$ are the \textit{null curvatures} of $\Sigma$. Instead, each leaf in the foliation constructed by Cederbaum and Sakovich in \cite{cederbaum} satisfies the equation 
\begin{equation*}
\textup{const}\equiv \mathcal H=\sqrt{H^2-P^2}=\sqrt{H+P}\sqrt{H-P}=:\sqrt{\Theta_+}\sqrt{\Theta_-}.
\end{equation*}
The two equations are different if the right hand side is strictly positive, as in the case of our interest. On the other hand, it is interesting to notice that they coincide if the right-hand side is zero, in which case we recover a well-known class called \textit{trapped surfaces}, or MOTS, which has been studied by various authors (for example \cite{eichmairmots} or \cite{eichmairmots2}).\\
\indent The existence of the foliations described above has been obtained by different methods through the various papers. The original approach by Huisken and Yau considered the volume preserving mean curvature flow starting from an Euclidean sphere with sufficiently large radius, showed global existence of the flow and obtained the CMC surface as the limit for large times. This flow was no longer studied in the context of asymptotically flat manifolds, except for the work of Corvino and Wu \cite{corvinowu}. The foliations were constructed instead by strategies based on the implicit function theorem (\cite{ye},\cite{huang}), a continuity method (\cite{metzger}, \cite{nerz1}, \cite{cederbaum})  or a Lyapunov-Schmidt reduction \cite{eichmairkoerber}. The volume preserving mean curvature flow has been recently studied by Sinestrari and the author \cite{vpmcf} in the context of asymptotically flat spaces, in order to recover the CMC-foliation under the optimal hypothesis of Nerz. We also mention the recent work \cite{guilisun} where an area preserving mean curvature flow is analyzed in the context of the Schwarzschild space.\\
\indent In this paper we study a modification of the volume preserving mean curvature flow. Let us describe our setting in more detail. We consider an initial data set $(\M,\overline \g,\overline\K)$ which is $C_{\frac12+\delta}^2$-asymptotically flat and constrained by the densities $(\overline\mu,\overline\J)$. For a general but fixed power $q\geq 2$ we set $\mathcal H=\sqrt[q]{H^q-|P|^q}$, whenever this quantity is defined. Notice that this is the Minkowski $q$-length of the vector $(\vec H,P)$, where $\vec H$ is the vector mean curvature. The \textit{volume preserving spacetime mean curvature flow} (VPSTMCF) is a family of time dependent immersions $F:\Sigma\times [0,T]\to\M$, with $\Sigma$ a closed 2-surface, which evolves according to 
\begin{equation}\label{stflow1}
\frac{\partial F}{\partial t}(t,\cdot)=-\left[\mathcal H(t,\cdot)-\hbar(t)\right] \nu(t,\cdot),
\end{equation}
where $\hbar(t)$ is the integral average of $\mathcal H (t,\cdot)$. Observe that, for $q=2$, $\mathcal H$ is the spacetime mean curvature defined in \cite{cederbaum}. As initial data for the flow \eqref{stflow1}, we consider a well-centered (in the sense of Nerz \cite{nerz1}) CMC-surface. As in \cite{huiskenyau}, the evolution is parametrized by a non-physical time parameter and takes place in a fixed spacelike slice, but in the present paper it has a speed that takes into account the spacetime texture of the initial data set. We aim to prove long-time existence of this flow, together with a convergence result. See the statement of Theorem \ref{mainthm} for details.
\\
\indent The approach we employ here to study the flow partly follows the one we have used in \cite{vpmcf} in the spacelike case.
We introduce a class of round surfaces defined by integral inequalities and prove that the solution
of the flow belongs to this class for every positive time. As in the construction in \cite{cederbaum}, which is based on the continuity method, we build the results on the spacelike case \cite{nerz1,vpmcf} by considering as initial value for the flow \eqref{stflow1} a surface with constant
(classical) mean curvature, which satisfies better estimates than a Euclidean coordinate sphere.
In this way we can assume on our solution a stronger condition on the oscillation of $H$ than in \cite{vpmcf},
see equation \eqref{almostCMCness}, which allows us to extend the spectral analysis of the stability operator
of \cite{cederbaum,nerz1} without the remainder which occurs in \cite{vpmcf}. On the other hand, the estimates
of the integral norms of the curvatures along the flow become more complicate and contain
nontrivial additional terms due to the nonlinearity of the speed. However, we are able to show that
the behaviour of the speed function is sufficiently close to linear for large radius of our surfaces
in order to estimate the additional terms and prove that our roundness class is preserved.
\\
\indent We now state the main result of the paper.
\begin{thm}\label{mainthm} Let $(\M,\overline\g,\overline\K)$ be a constrained $C_{\frac12+\delta}^2$-asymptotically flat initial data set in the sense of Definition \ref{def1} and \ref{def2}, and suppose that $E_\textup{ADM}>0$. Fix $q\geq 2$. Let $\iota:\Sigma\hookrightarrow \M$ be a \textup{closed CMC-surface immersed} in $(\M,\overline \g)$ and, setting $\sigma=\sigma_\Sigma:=\sqrt{|\Sigma|/4\pi}$, suppose that there exists $C_0>0$ such that  
\begin{equation}\label{hypstmainth2in}
\|\overset{\circ}{A}\|_{L^4(\Sigma)}\leq C_0\sigma^{-1-\delta}, \qquad |\vec z_\Sigma|\leq C_0\sigma^{1-\delta}, \qquad \frac{\sigma}{r_\Sigma}\leq 1+C_0^{-1},
\end{equation}
where $r_\Sigma:=\displaystyle\min_{x\in\Sigma}|\vec x(x)|$. Then, there exists $\sigma_0=\sigma_0(C_0,\overline c,\delta,q)>1$ such that if $\sigma>\sigma_0$, the solution $\Sigma_t$ to the spacetime mean curvature flow \eqref{stflow1} exists for every $t\in[0,\infty)$ and converges exponentially fast to a surface $\Sigma^{\textup{st}}_\infty$ satisfying the prescribed mean curvature equation
\begin{equation}\label{eqqstinfty}
H_{\Sigma_\infty^{\textup{st}}}^q=|P|_{\Sigma_\infty^{\textup{st}}}^q+\hbar_{\Sigma_\infty^{\textup{st}}}^q
\end{equation}
for some constant $\hbar_{\Sigma_\infty^{\textup{st}}}>0$.
\end{thm}
This flow approach provides an alternative construction of the foliation obtained in \cite{cederbaum}. In particular, when $q=2$, letting evolve the CMC-foliation $\{\Sigma^\sigma\}_{\sigma\geq \sigma_0}$ constructed in \cite{nerz1} or in \cite{vpmcf}, we obtain when $t\to+\infty$ a CSTMC-foliation $\{\Sigma_{\textup{st}}^\sigma\}_{\sigma\geq\sigma_0}$ which coincides with the one introduced by Cederbaum-Sakovich. Our construction has an independent interest in the analysis of the behaviour of curvature flows in asymptotically flat spaces and proves that the foliation exists also when $q>2$. While in the case $q=2$ equation \eqref{eqqstinfty} is related to the conjugate momentum tensor $\overline\pi:=\left(\textup{tr}_{\overline\g}\overline\K\right)\overline\g-\overline\K$,
see \cite[Lemma 10]{cederbaum}, no similar relation is known, to the knowledge of the author, for general $q$. Thus the physical interpretation of the center of mass of the CSTMC-foliation, when existing, is not clear, at least when $q$ is greater than but close to 2.
\\
\indent We finally recall that the Ricci flow plays an important role in the mathematics of high energy theoretical physics, see \cite{carfora}. In our opinion, extrinsic curvature flows like the ones considered in this paper coupled with a Ricci flow evolution of a background metric may prove to be useful in physical applications. A further motivation for our study comes from the recent work of Huisken and Wolff, where the authors studied a spacetime version of the inverse mean curvature flow, see \cite{huiskenwolff}. Finally, an area preserving spacetime mean curvature flow, in the context of the Schwarzschildean lightcones, has been studied by Kr\"oncke and Wolff in \cite{krwo}.
\section*{Acknowledgements}
\ \ \ \ \ \ \ The author acknowledges support from the MIUR Excellence Department Project, CUP E83C18000100006, and the MUR Excellence Department Project MatMod@TOV, CUP E83C23000330006, awarded to the Department of Mathematics, University of Rome Tor Vergata, and from the MUR Prin 2022 Project "\textit{Contemporary perspectives on geometry and gravity}" CUP E53D23005750006. The author would like to thank his Ph.D. advisor Carlo Sinestrari for suggesting the topic, his constant supervision and the helpful comments during the writing of this paper.
\section{Preliminaries and notations}
\subsection{Definitions and basic properties} A \textit{spacetime}, say $\mathcal M$, is a connected 4-manifold equipped with a $(0,2)$-type tensor field $\textbf{\textup{g}}$ non-degenerate and \textit{Lorentzian}, i.e. with signature $(-1,1,1,1)$, which satisfies the so called \textit{Einstein} (constraint) \textit{equations}
\begin{equation}\label{einsteineq}
\textbf{\textup{Ric}}-\left(\frac{\textbf{\textup{S}}}{2}\right)\textbf{\textup{g}}=\textbf{\textup{T}},
\end{equation}
where $\textbf{\textup{Ric}}=\textbf{\textup{Ric}}(\cdot,\cdot)$ is the Ricci tensor of $(\mathcal{M},\textbf{\textup{g}})$, $\textbf{\textup{S}}=\textup{tr}_{\textbf{\textup g}}\textbf{\textup{Ric}}$ and $\textbf{\textup{T}}$ is a symmetric smooth $(0,2)$-type tensor field called \textit{energy-momentum tensor}. In the following, we will be interested in 3-dimensional submanifolds $\M$ of a spacetime such that the restriction of the Lorentzian metric $\textbf g$ to $\M$ is a Riemannian metric. These submanifolds are known as \textit{spacelike submanifolds}. We will always require $(\mathcal M, \textbf g)$ to be \textit{time-orientable}. In particular, we suppose that there exists a globally defined timelike smooth vector field $\boldsymbol e_0$. We also denote by $\boldsymbol A$ the second fundamental form of the immersion $\boldsymbol j:\M\hookrightarrow \mathcal M$. We then define the metric $\overline \g:=\boldsymbol j^*\textbf g$ and 
\begin{equation}\label{initialdataconstandk}
\overline\K(\cdot,\cdot):=\langle \boldsymbol A(\cdot,\cdot),-\boldsymbol e_0\rangle, \qquad \overline\mu:=\textbf T(\boldsymbol e_0,\boldsymbol e_0), \qquad \overline\J(\cdot):=\textbf T(\boldsymbol e_0,\cdot),
\end{equation}
which are, respectively, the (scalar) \textit{spacetime second fundamental form} of $\M$, the \textit{energy density} and the \textit{momentum density}. Through the Gauss-Mainardi-Codazzi equations, it can be proved that \eqref{einsteineq} implies the the system \eqref{einstconstsyst} holds. Conversely, Choquet-Bruhat \cite{bruhat} proved, in the vacuum and some special cases, that the validity of the system \eqref{einstconstsyst} implies the existence of a spacetime $(\mathcal M,\textbf g)$ such that $(\M,\overline\g,\overline\K)$ sits into $(\mathcal M,\textbf g)$ and \eqref{initialdataconstandk} are satisfied. See also \cite{glockle} for relations with the (strict) dominant energy condition. This remark allows us to define an initial data set as a notion independent from the one of spacetime manifold (at least formally), but which includes the whole information that a spacetime encodes.
\\
\indent 
In the following, we will always assume that $(\M,\overline\g,\overline\K)$ is an asymptotically flat initial data set in the sense of Definition \ref{def1}. We denote by $\overline\Gamma_{\alpha\beta}^\gamma$ the  Christoffel symbols of the Riemannian manifold $(\M,\overline\g)$, and by $\overline\nabla$ the Riemannian connection. By $\overline{\textup{Rm}}$, $\overline{\textup{Ric}}$ and $\overline{\textup S}$ we denote the Riemann and Ricci curvature tensors on $(\M,\overline \g)$. We denote the Euclidean coordinate spheres on $(\M,\overline \g)$ by $\mathbb S_R(\vec z)$ , for some $R>1$ and $\vec z\in\R^3$, by an abuse of notation using that $\M\setminus \textup C$ is diffeomorphic to $\R^3\setminus \overline{\mathbb B}_1(\vec 0)$.\\
\indent We say that $(\M,\overline \g)$ satisfies the \textit{mass condition} if the scalar curvature satisfies the inequality $|\overline{\textup S}|\leq \overline c|\vec x|^{-3-\delta}$, possibly enlarging the constant $\overline c$. We remark that, in the case of constrained initial data sets, the mass condition is implied by \eqref{einstconstsyst} together with \eqref{constrainedmuJ}. This is the only reason we need \textit{constrained} initial data sets. We will use the following equivalent characterization of the ADM-energy introduced by Arnowitt, Deser and Misner in \cite{adm}. The equivalence of the definitions has been proved by Miao and Tam in \cite{miaotam}.
\begin{dfn}\label{adm-energy} Let $(\M,\overline\g)$ be a $C_{\frac12+\delta}^2$-asymptotically flat 3-manifold that satisfies the \textit{mass condition}. The ADM-\textit{energy} is defined as 
\begin{equation}\label{admenergy}
E_{\textup{ADM}}:=-\lim_{R\to\infty}\frac{R}{8\pi}\int_{\vec x^{-1}(\mathbb S_R(\vec 0))}\overline{\textup{G}}(\nu_R,\nu_R) \ d\mu_R,
\end{equation}
where $\overline{\textup{G}}:=\overline{\textup{Ric}}-\left(\frac{\overline{\textup{S}}}{2}\right)\overline \g$ is the (spacelike) Einstein tensor, $\vec x^{-1}(\mathbb{S}_R(\vec 0))$ is the Euclidean sphere immersed in $(\M,\overline\g)$, and $\nu_R$ and $d\mu_R$ are, respectively, its normal vector and its volume form in $(\M,\overline \g)$.
\end{dfn}
In the following, by \textit{surface} we mean an immersion $\iota:\Sigma\hookrightarrow \M\setminus \textup C$, with $\dim\Sigma=2$, which is closed, connected and 2-faced. Since $\M\setminus \textup C$ is diffeomorphic to $\R^3\setminus\overline{\mathbb B}_1(\vec 0)$, the surface $\Sigma$ inherits two Riemannian metrics: a \textit{physical} metric $g:=\iota^*\overline\g$ and a \textit{Euclidean} metric $g^e:=\iota^*\overline \g^e$, where $\overline \g^e$ is the Euclidean metric on $\M$. From now on, we will use the apex $e$ each time a quantity is computed with respect to the Euclidean metric, and we will omit the apex if it is computed using the physical metric. Then, fixed an outer unit normal $\nu:\Sigma\to T\M$, we represent with $A=\{h_{ij}\}$, $H$ and $d\mu$, respectively, the second fundamental form, the mean curvature and the volume form of $\Sigma$ with respect to $g$. On the other hand, if $\nu^e:\Sigma\to T\M$ is the outer unit normal of $\Sigma$ with respect to $\overline\g^e$, we represent the same quantities with $A^e=\{h^e_{ij}\}$, $H^e$ and $d\mu^e$. Observe that, when we are on a hypsersurface $\Sigma$, we use the latin indexes $i,j,k,l$, etc, to distinguish from the ambiental coordinates, which are indicated with the greek indexes $\alpha,\beta,\gamma,\epsilon$, etc. Finally, we define 
\begin{equation*}
h:=\frac1{|\Sigma|}\int_\Sigma H\ d\mu,\qquad h^e:=\frac1{|\Sigma^e|}\int_{\Sigma} H^e \ d\mu^e,
\end{equation*}
which are, respectively, the mean of the mean curvature computed with respect to the physical and the Euclidean metric. Here $|\Sigma|=\int_\Sigma d\mu$ and $|\Sigma^e|=\int_\Sigma d\mu^e$. The asymptotic flatness of the ambient manifold allows to deduce the following estimates, which can be proved by standard computations, see e.g. \cite[Lemma 1.5]{metzger}, \cite[Lemma 11]{cederbaum}.
\begin{lem} Let $\iota:\Sigma \hookrightarrow \M$ be a surface immersed in a $C_{\frac12+\delta}^2$-flat 3-manifold $\M$. Then there exists $C>0$, only depending on the constant $\overline c$ in Definition \ref{def1}, such that
\begin{equation}
|g-g^e|\leq C|\vec x|^{-\frac12-\delta}, \qquad |\Gamma_{ij}^k-\left(\Gamma^e\right)_{ij}^k|\leq C|\vec x|^{-\frac32-\delta},
\end{equation}
\begin{equation}
|d\mu-d\mu^e|\leq C|\vec x|^{-\frac12-\delta}d\mu,
\end{equation}
\begin{equation}
|\nu-\nu^e|\leq C|\vec x|^{-\frac12-\delta}, \qquad |\nabla\nu-\nabla^e\nu^e|\leq C|\vec x|^{-\frac32-\delta},
\end{equation}
\begin{equation}\label{AAeineq}
|A-A^e|\leq C\left(|\vec x|^{-\frac32-\delta}+|\vec x|^{-\frac12-\delta}|A^e|\right),
\end{equation}
\begin{equation}
|\nabla A-\nabla^eA^e|\leq C\left(|\vec x|^{-\frac52-\delta}+|\vec x|^{-\frac12-\delta}|\nabla^eA^e|\right).
\end{equation}
In addition, if $|A|\leq 10|\vec x|^{-1}$, then 
\begin{equation}
|H-H^e|\leq C|\vec x|^{-\frac32-\delta}, \qquad |\overset{\circ}{A}-\overset{\circ}{A^e}|\leq C|\vec x|^{-\frac32-\delta}.
\end{equation}
\end{lem}
We also define the barycenter of $\Sigma$ as 
\begin{equation*}
\vec z_\Sigma:=\int_\Sigma \iota \ d\mu.
\end{equation*}
Finally, in order to estimate the Euclidean position of an immersed surface and its area, we introduce the following definitions.
\begin{dfn}\label{defradii} Let $(\M,\overline\g)$ be a 3-manifold, and consider an immersed surface $\iota:\Sigma\to \M\setminus \textup C$ with induced metric $g=\iota^*\overline\g$. Then we set 
\begin{equation*}
r_\Sigma:=\min_{x\in\Sigma}|\vec x(\iota(x))|, \qquad R_\Sigma:=\max_{x\in\Sigma}|\vec x(\iota(x))|, \qquad \sigma_\Sigma:=\sqrt{\frac{|\Sigma|_g}{4\pi}}.
\end{equation*}
These \textit{radii} are called \textit{Euclidean radius}, \textit{Euclidean diameter} and \textit{area radius}, respectively. 
\end{dfn}
As in the previous literature, see \cite{nerz1}, \cite{cederbaum}, we define the Sobolev norms on $\Sigma$ as follows 
\begin{equation*}
\|f\|_{W^{0,p}(\Sigma)}:=\|f\|_{L^p(\Sigma)}, \qquad \|f\|_{W^{k+1,p}(\Sigma)}:=\|f\|_{L^p(\Sigma)}+\sigma_\Sigma \|\nabla f\|_{W^{k,p}(\Sigma)}
\end{equation*}
for $p\in [1,\infty]$ and $k\in\N^{\geq 0}$. As usual, we set $H^k=W^{2,k}$.
\subsection{Surfaces in initial data sets} If $(\M,\overline\g,\overline\K)$ is an asymptotically flat initial data set and $\Sigma\hookrightarrow \M$, the surface $\Sigma$ can be seen as a 2-codimensional submanifold of the 4-dimensional spacetime manifold associated to $(\M,\overline\g,\overline\K)$. Since in this case we have the decomposition of the tangent space of a point $p\in\mathcal M$ given by $T_p\mathcal M=\langle \boldsymbol e_0\rangle_p\oplus^{\perp_\mathcal M} T_p\Sigma\oplus^{\perp_\M}N_p \Sigma$. This decomposition allows to define the so called \textit{null mean curvatures}, which are given by
\begin{equation}
\theta^\pm=P\pm H,
\end{equation}
where $P:=\textup{tr}_g(\overline\K)\equiv g^{ij}\overline\K_{ij}$ is the so called \textit{timelike mean curvature}. At this point, we define the \textit{spacetime mean curvature}, if it exists, as 
\begin{equation*}
\mathcal{H}:=\sqrt{H^2-P^2}.
\end{equation*}
It is essentially the Minkowski 2-length of the vector $(\vec H,P)$, where $\vec H$ is the \textit{vector mean curvature} of $\Sigma$. In this paper, we will mainly consider the Minkowski $q$-norm $|(\vec H,P)|_q=\left(H^q-|P|^q\right)^\frac1q$, with $q\geq 2$. Whenever there are no ambiguities, for a fixed $q\geq 2$ we will continue to call this quantity $\mathcal H$. In all the cases, we set 
\begin{equation*}
\hbar:=\fint_\Sigma \mathcal{H} \ d\mu.
\end{equation*}
In the following preliminary Lemma, we show some properties of surfaces in initial data sets.
\begin{lem}\label{propertiesinitialdata} Let $(\M,\overline\g,\overline \K)$ be a $C_{\frac12+\delta}^2$-asymptotically flat initial data set. Then there exist constants $C=C(\overline c)>0$ and $\cin=\cin(\overline c)>0$, also depending on the choice of $q$, such that if $\iota:\Sigma\hookrightarrow \M\setminus \textup C$ is a \textup{surface} with induced metric $g:=\iota^*\overline\g$ and there exists $\sigma>1$ such that
\begin{equation}\label{eq222}
2\sigma\geq r_\Sigma\geq \frac{\sigma}2, \qquad \frac{1}{\sigma}\leq H_x\leq \frac{\sqrt{5}}{\sigma} \quad \forall x\in\Sigma,
\end{equation}
then $\mathcal H_x$ is well defined for every $x\in \Sigma$ and the following inequalities hold:
\begin{equation}\label{eq222bis}
\|P\|_{L^\infty(\Sigma)}+\sigma\|\nabla P\|_{L^\infty(\Sigma)}\leq C\sigma^{-\frac32-\delta},
\end{equation}
\begin{equation*}
\sup_\Sigma |\mathcal H-H|\leq C\sigma^{-1-\frac12q-q\delta}, \qquad |h-\hbar|\leq C\sigma^{-1-\frac12q-q\delta}.
\end{equation*}
Finally, if $H$ is constant on $\Sigma$, i.e. $H\equiv h$, then $\|\mathcal H-\hbar\|_{L^\infty(\Sigma)}+\sigma\|\nabla\mathcal H\|_{L^\infty(\Sigma)}\leq \cin\sigma^{-1-\frac12q-q\delta}$.
\end{lem}
In the following, we will mainly use the $H^1$-estimate on $\mathcal H-\hbar$, which follows from the $W^{1,\infty}$ bound in the above statement, and we will continue to call $\cin$ the constant at the right-hand side. Observe that this is the only case in which we use the lowercase in order to indicate a constant depending on the setting and not on the "roundness of the surface" (in a sense we will make more clear later, see Definition \ref{roundsurface} below).
\begin{proof}
Equation \eqref{eq222bis} follows from the fact that $|P_x|\leq 2|\overline{\textup{K}}|_{\overline \g}\leq 2\overline c r_\Sigma^{-\frac32-\delta}=\textup O(\sigma^{-\frac32-\delta})$,
\begin{equation*}
(\nabla P)_i=g^{ik}g^{jl}\nabla_k\overline{\textup{K}}_{jl},
\end{equation*}
and the fact that $|\nabla_k\overline{\textup{K}}_{jl}|\leq \overline c\sigma^{-\frac52-\delta}$. The other inequalities follow from the Lagrange mean value theorem and the identity
\begin{equation}\label{eq18b}
\mathcal{H}^{q-1}\nabla\mathcal{H}=H^{q-1}\nabla H-|P|^{q-1}\left(\frac{P}{|P|}\right)\nabla P.
\end{equation}
\end{proof}
\subsection{Round surfaces} 
We now introduce a class of surfaces, which are close to a Euclidean sphere of the same area radius in a precise (quantitative) way. The aim is to find a class which
is invariant under the volume preserving spacetime mean curvature flow under an appropriate choice of the
parameters and for large enough radius. Our class of round surfaces coincide with the one introduced in \cite{vpmcf}. Other classes of round surfaces, which are related to
the methods used there, have been introduced in \cite{huiskenyau}, \cite{metzger}, \cite{nerz1}, \cite{cederbaum}.
\begin{dfn}\label{roundsurface} Let $(\M,\overline \g,\overline{\textup{K}})$ be a $C_{\frac12+\delta}^2$-asymptotically flat initial data set and let $\iota:\Sigma\to M$ be a surface.\\
\noindent (i) For a given approximate radius $\sigma>1$ and parameters $\eta>0$, $B_1$, $B_2>0$, we say that $(\Sigma,g)$ is a \textit{round surface} in $(\M,\overline\g,\overline \K)$, and we write quantitatively $\Sigma\in \mathcal{W}_\sigma^\eta(B_1,B_2)$ if the following inequalities are satisfied
\begin{equation}\label{radius1}
|A|<\sqrt{\frac{5}{2\sigma^2}}, \qquad \kappa_i\geq \frac1{2\sigma},
\end{equation}
\begin{equation}\label{radius2}
\frac{7}2\pi\sigma^2<|\Sigma|_g<5\pi\sigma^2, \qquad \frac34\leq \frac{r_\Sigma}{\sigma}\leq \frac{R_\Sigma}{\sigma}\leq \frac54,
\end{equation}
\begin{equation}\label{cond2defroundst}
\|\overset{\circ}{A}\|_{L^4(\Sigma,\mu)}<B_1\sigma^{-1-\delta},
\end{equation}
\begin{equation}\label{cond3defroundst}
\eta\sigma^{-4}\|\mathcal H-\hbar\|_{L^4(\Sigma)}^4+\|\nabla \mathcal H\|_{L^4(\Sigma)}^4<B_2\sigma^{-8-4\delta}.
\end{equation} 
(ii) For given $\sigma>1$ and $\eta$, $B_1$, $B_2$, $B_\textup{cen}>0$, we say that $(\Sigma,g)$ is a \textit{well-centered round surface}, and we write $\Sigma\in{\mathcal{B}}_\sigma^\eta(B_1,B_2,B_\textup{cen})$ if it satisfies the above properties and in addition 
\begin{equation}
|\vec z_\Sigma|<B_\textup{cen}\sigma^{1-\delta}
\end{equation}
\end{dfn}
The decay rates in conditions \eqref{cond2defroundst}-\eqref{cond3defroundst} are modelled on the properties of the leaves of the CMC-foliation constructed by Nerz in \cite{nerz1}. See also the construction proposed in \cite{vpmcf}. In particular, in his Theorem 5.1, Nerz proved the existence of an exhaustive family of constant mean curvature surfaces which foliate an asymptotically flat manifold with non-zero ADM-energy. We denote by $\{\Sigma^s\}_{s\geq s_0}$, for a certain $s_0>1$, the CMC-foliation constructed by Nerz in \cite{nerz1}. We avoid to use $\sigma$ as a parameter of this foliation as Nerz does since in the present paper it plays a different role, as we can see in Definition \ref{roundsurface}. Nerz proved that this foliation satisfies 
\begin{equation}\label{nerzconds}
H^{\Sigma^s}=\frac2s, \qquad \|\overset{\circ}{A}\|_{H^1(\Sigma^s)}\leq C_\textup{Nerz} s^{-\frac32-\delta}, \qquad |\vec z_{\Sigma^s}|\leq C_\textup{Nerz} s^{1-\delta},
\end{equation}
for some $C_\textup{Nerz}>0$. Moreover, \cite[Prop. 4.4]{nerz1} proves that $|s-\sigma_s|\leq C\sigma_s^{\frac12-\delta}$, where $\sigma_s:=\sigma_{\Sigma^s}$. Note also that \cite[Prop. 4.4]{nerz1}, combined with \eqref{nerzconds}, implies
\begin{equation}
\sigma_{\Sigma^s}-C\sigma_{\Sigma^s}^{1-\delta}\leq |\vec x|=|\vec z^s+\sigma_{\Sigma^s}\nu^s+f^s\nu^s|\leq \sigma_{\Sigma^s}+C\sigma_{\Sigma^s}^{1-\delta},
\end{equation}
that is $|r_{\Sigma^s}-\sigma_{\Sigma^s}|\leq C\sigma_{\Sigma^s}^{1-\delta}$. Then 
\begin{equation}\label{remark26}
\frac{r_{\Sigma^s}}{\sigma_s}\geq 1-C\sigma_s^{-\delta}.
\end{equation}
Thus, for $s$ large, this foliation satisfies \eqref{hypstmainth2in} with $\sigma=\sigma_s$. Finally observe that, by Lemma \ref{propertiesinitialdata}, the leaves $\Sigma^s$ also satisfy
\begin{equation*}
\|\mathcal H-\hbar\|_{W^{1,2}(\Sigma^s)}\leq \cin\sigma^{-\frac{q}2-q\delta},
\end{equation*}
for some $\cin=\cin(\overline c)>0$. For this reason, we will use a fixed leaf of Nerz's foliation as the initial datum of our flow. \\
\\
\noindent \textbf{Notation for the constants.} Throughout the paper, when deriving estimates on geometric quantities on a surface $\Sigma$, we denote by $C,C_1,C_2,...$ constants which only depend on properties of the ambient manifold, such as $\overline c,\delta$ in Definition \ref{def1} or the energy $E_\textup{ADM}$ and by $c,c_1,c_2,\tilde c$,... constants which in addition depend on the constants $B_1,B_2,B_\textup{cen}$ in the previous conditions. We say that a constant is universal if it is independent on any other parameter of our problem. As
usual, the letters $c$ or $C$ will often denote constants which may change from one line to the other, but each time depending on the same parameters.
\begin{rem}\label{preliminaryremark27}
Property \eqref{radius2} implies that the three radii of Definition \ref{defradii} are comparable among each other and with $\sigma$.
In particular this property implies, because of the asymptotic flatness of $\M$ as in \eqref{asymptc8}, the bound on the Riemann tensor
\begin{equation}\label{boundonriem}
 |\overline{\textup{Rm}}| \leq C(\bar c) \sigma^{-\frac 52 -\delta} \mbox{ on } \Sigma.
\end{equation}
\end{rem}
\indent In the next Lemma, we collect some well-known properties of round surfaces. See \cite[Remark 2.6]{vpmcf} and \cite[Lemma 2.7]{vpmcf} for more details and a proof.
\begin{lem}\label{cor1} Let $(\M,\overline \g)$ be a $C_{\frac12+\delta}^2$-asymptotically flat manifold. Fix any $\eta$, $B_1$ and $B_2>0$. Then there exists $\sigma_0=\sigma_0(B_1,B_2,\eta,\overline c,\delta)>0$ such that any surface which belongs to $\mathcal{W}_\sigma^\eta(B_1,B_2)$ for some $\sigma>\sigma_0$ satisfies the following properties:
\begin{enumerate}[label=\textup{(\roman*)}]
\item There exists $\cs>0$ such that it holds 
\begin{equation}\label{sobolevineq}
\|\psi\|_{L^2}\leq\frac{\cs}{\sigma}\|\psi\|_{W^{1,1}(\Sigma)}, \qquad \forall \ \psi\in W^{1,1}(\Sigma)
\end{equation}
and, for every $p>2$,
\begin{equation}\label{othersob11}
\|\psi\|_{L^\infty}\leq 2^{\frac{2(p-1)}{p-2}}\cs\sigma^{-\frac{2}{p}}\|\psi\|_{W^{1,p}}, \qquad \forall \psi\in W^{1,p}(\Sigma)
\end{equation}
Moreover, there exists a constant $c(B_2,\eta)>0$ such that  
\begin{equation}\label{LinfcontrolH}
\|H-h\|_{L^\infty}\leq c(B_2,\eta)\sigma^{-\frac32-\delta}.
\end{equation}
\item It holds the estimate 
\begin{equation}\label{eqbis529}
\left|h-\frac2\sigma_\Sigma\right|\leq c(B_1,B_2,\overline c)\sigma^{-\frac32-\delta}.
\end{equation}
\item There exists a constant $B_\infty=B_\infty(B_1,B_2,\eta,\delta,\overline c)$ such that $\|\overset{\circ}{A}\|_{L^\infty(\Sigma)}\leq B_\infty\sigma^{-\frac32-\delta}$.
\item There exists $c=c(\delta,\overline c,B_1,B_2,\eta)$, $c_0=c(B_1,\overline c,\delta)$, $\vec z_0\in\mathbb R^3$, and $f:\mathbb{S}_{\sigma_\Sigma}(\vec z_0)\to \mathbb{R}$ such that 
\begin{equation}\label{eq2110}
\Sigma^e=\textup{graph}(f), \qquad \|f\|_{W^{2,\infty}}\leq c\sigma^{\frac12-\delta}, \qquad |\vec z_0-\vec z_\Sigma|\leq c_0\sigma^{\frac12-\delta}.
\end{equation}
\item It holds $\|A\|_{L^4(\Sigma)}\leq 4^\frac14\sqrt{5/2}\pi^\frac14 \sigma^{-\frac12}$ and thus there exists a constant $\cper=\cper(\delta,\overline c)$ such that 
\begin{equation*}
\|H-h\|_{L^4(\Sigma)}\leq \cper\|\overset{\circ}{A}\|_{L^4(\Sigma)}+\cper\sigma^{-1-\delta}.
\end{equation*}
\end{enumerate}
\end{lem}
\begin{rem}\label{remark28} Most of the results in the previous Lemma follows from the fundamental result of DeLellis-Müller \cite{delellismuller1}. We also remark that, by Theorem 1.1 in \cite{delellismuller1}, the assumptions in \eqref{hypstmainth2in} imply that $\frac{\sigma}{r_\Sigma}$ is also bounded away from zero. In fact, combining \eqref{AAeineq} with the latter assumption in \eqref{hypstmainth2in}, and using the DeLellis-Müller's Theorem, Point (iv) of Lemma \ref{cor1}, combined with the second assumption in \eqref{hypstmainth2in}, implies that $\frac{\sigma}{r_\Sigma}\geq 1-C_0^{-1}$, provided that $\sigma$ is suitably large.
\end{rem}
\section{Spectral theory}\label{spthids}
\subsection{Mass and stability operator} In this section we consider closed surfaces $\Sigma$ belonging to a roundness
class ${\mathcal W}_\sigma^\eta(B_1,B_2)$ for fixed parameters $\eta$, $B_1$, $B_2$ and a general large $\sigma$. Moreover, we will suppose that $E_\textup{ADM}>0$ and that $\Sigma$ satisfies 
\begin{equation}\label{almostCMCness}
\|H-h\|_{L^2(\Sigma)}\leq \cin\sigma^{-1-2\delta}.
\end{equation}
We sometimes say that a surface $\Sigma$ satisfying \eqref{almostCMCness} is $(\cin,\sigma)$-\textit{almost CMC}. We simply say that $\Sigma$ is $\cin$-\textit{almost CMC} if $\Sigma$ is $(\cin,\sigma_\Sigma)$-almost CMC. Thus, we will tacitly mean that the constants $c$ and $\sigma_0$ which appear in the statements below only depend on $\eta,B_1,B_2$, on the constants $\cin,\overline c,\delta$ and possibly on the energy $E_\textup{ADM}$.\\
\indent Since our aim is mostly to investigate the extrinsic geometry of surfaces, we recall a (quasi-) local notion of mass, which goes back to \cite{hawking68}. 
\begin{dfn}[Hawking energy] Let $(\M,\overline \g)$ be a 3-dimensional manifold, and $\iota:\Sigma\hookrightarrow \M$ be a surface. Let $g:=\iota^*\overline \g$ be the induced metric. The \textit{Hawking energy} of $\Sigma$ is defined as
\begin{equation}\label{hawmass311}
m_H(\Sigma):=\sqrt{\frac{|\Sigma|_g}{16\pi}}\left(1-\frac{1}{16\pi}\int_\Sigma H^2 \ d\mu\right).
\end{equation}
\end{dfn}
\begin{rem} The original definition the Hawking energy involves the spacetime mean curvature $\mathcal H$ instead of the (spacelike) mean curvature $H$. Precisely, the definition in \eqref{hawmass311} takes the name of \textit{Geroch mass}. However, it is well-known that the two definition are asymptotic for large round surfaces (see \cite{cederbaum}) and thus we continue to use the (local) notion of mass in \eqref{hawmass311}.
\end{rem}
If $\Sigma\in \mathcal W_\sigma(B_1,B_2)$, the Gauss-Bonnet theorem shows that the Hawking mass is asymptotic to the ADM-energy when considering round surfaces of large radius, proving in particular the following Lemma. See \cite[Appendix A]{nerz1} for a proof.
\begin{lem} There exist $\tilde c$ and $\sigma_0$ such that, for every $\Sigma\in\mathcal W_\sigma^\eta(B_1,B_2)$ with $\sigma\geq \sigma_0$ we have
\begin{equation}\label{mHE43}
\left|m_{H}(\Sigma)+\frac{\sigma_\Sigma}{8\pi}\int_\Sigma\overline{\textup{G}}(\nu,\nu) \ d\mu\right|\leq \tilde c\sigma^{-\delta}, \qquad \left|E_{\textup{ADM}}+\frac{\sigma_\Sigma}{8\pi}\int_\Sigma \overline{\textup{G}}(\nu,\nu) \ d\mu\right|\leq \tilde c\sigma^{-\delta}.
\end{equation}
\end{lem}
We now introduce the stability operator, which occurs as the second variation of the area
functional.
\begin{dfn}
Given a surface $\iota:\Sigma\hookrightarrow \M$ and a smooth function $\textup f\in H^2(\Sigma)$, we define the \textit{stability operator} associated to $\Sigma$, $L^\Sigma:H^2(\Sigma)\to L^2(\Sigma)$, as 
\begin{equation}\label{stabilityop}
L^\Sigma\textup f:=-\Delta \textup f-(|A|^2+\overline{\textup{Ric}}(\nu,\nu))\textup f.
\end{equation}
\end{dfn}
We simply write $L$ instead of $L^\Sigma$ whenever the role of the surface $\Sigma$ is not ambiguous. In \cite{nerz1} and \cite{cederbaum} the spectral properties of $L$ are studied in the case of CMC-surfaces. In particular, they show that $L$ is invertible if $E_\textup{ADM}\neq 0$ and positive definite on functions with zero mean if $E_\textup{ADM}>0$. Here we generalize this analysis to round surfaces where we only assume that $H$ has a small oscillations. See \cite[Section 3]{vpmcf} for the general case of round surfaces. In the context of the initial data sets, Cederbaum-Sakovich \cite{cederbaum} proved results analogous to those of Nerz, modifying the operator $L^\Sigma$ with additional terms. Even if we work in initial data sets, we continue to use the definition of stability operator given in \eqref{stabilityop}. The reason is that, as we will see in the next Sections, we will obtain the stability operator in computing the evolution of some geometric quantities. In estimating the terms involved in these computations, the main addend will be the operator in \eqref{stabilityop}, while the remainder will be small with a precise order of decay (for round surfaces with large radius).\\
\indent Since the "differential part" of the stability operator is totally given by the Laplace-Beltrami operator, in order to understand the properties of $L$ we are interested in we have to briefly review the spectral theory for the operator $-\Delta$. The following Lemma is taken from \cite[Lemma 2]{cederbaum}. Observe that, together with adapting the notations of the Lemma with our definition of roundness class, we remove the hypothesis of having a CMC-surface. In fact, reading the proof of \cite[Lemma 2]{cederbaum} with attention, one can observe that the CMC-hypothesis is just needed in order to compare the area radius with the \textit{curvature radius} used in \cite{nerz1} and \cite{cederbaum}.
\begin{rem} At the light of Lemma \ref{cor1}, i.e. of the De Lellis-Müller theorem \cite[Thm. 1.1]{delellismuller1}, scalar functions on a round surface $\Sigma$ can be also meant as functions on the approximating sphere $\mathbb{S}_{\sigma_\Sigma}$. With an abuse of notation, we identify such kind of functions.
\end{rem}
We first recall some properties of the Laplace-Beltrami operator on a round sphere $\mathbb S_\sigma(\vec 0) \subset \R^3$ with the Euclidean metric. 
On a general closed surface, the eigenvalues of the Laplace operator are all positive, except the first one which is zero, with eigenspace given by the constant functions. For the Euclidean sphere, the first nonzero eigenvalue has multiplicity three and is given by
$$
\lambda^e_\alpha=\frac{2}{\sigma^2}, \ \alpha=1,2,3.
$$
An orthonormal basis for the eigenspace is given by the normalized coordinate functions 
$$
f_\alpha^e(\vec x)  =\sqrt{\frac{3}{4\pi \sigma_\Sigma^4}} \vec x_\alpha, \ \alpha=1,2,3,
$$
restricted on $\mathbb S_\sigma(\vec 0)$. The remaining eigenvalues satisfy the bound
$$
\lambda^e_i \geq \lambda^e_4=\frac{6}{\sigma^2}, \ \forall \, i \geq 4.
$$
We recall the statement of Lemma 2 of \cite{cederbaum}, which measures how much the first eigenvalues and the corresponding eigenfunctions of the Laplace-Beltrami operator on a round surface in the physical metric differ from the ones of the approximating sphere in the Euclidean metric.
\begin{lem}\label{ortonormalsystemL2} There exist $c>0$ and $\sigma_0>1$ such that, if $\Sigma\in {\mathcal W}_\sigma^\eta(B_1,B_2)$ with $\sigma\geq \sigma_0$, there is a complete orthonormal system in $L^2(\Sigma)$ consisting of the eigenfunctions $\{f_\alpha\}_{\alpha=0}^\infty$ such that 
\begin{equation*}
-\Delta f_\alpha=\lambda_\alpha f_\alpha, \qquad \text{with $0=\lambda_0<\lambda_1\leq \lambda_2\leq...$}
\end{equation*}
Set $\mathbb{S}{\sigma_\Sigma}$ to be the round sphere approximating $\Sigma$ in the sense of \textup{Lemma \ref{cor1}}. Then there exists an orthonormal triple $\{f_1^e,f_2^e,f_3^e\}$ of eigenfunctions of $-\Delta^{\mathbb S_{\sigma_\Sigma}}$ such that, for $\alpha=1,2,3$,  
\begin{equation*}
\left|\lambda_\alpha-\frac{2}{\sigma_\Sigma^2}\right|\leq c\sigma^{-\frac52-\delta}, \quad \|f_\alpha-f_\alpha^e\|_{W^{2,2}(\Sigma)}\leq c\sigma^{-\frac12-\delta}, 
\quad \left\|\overset{\circ}{\textup{Hess}}(f_\alpha)\right\|_{L^2(\Sigma)}\leq  c\sigma^{-\frac52-\delta},
\end{equation*}
Moreover 
\begin{equation}\label{ineq12stim}
\int_\Sigma \left|\langle \nabla f_\alpha,\nabla f_\beta\rangle-\frac{3\delta_{\alpha\beta}}{\sigma_\Sigma^2|\Sigma|_g}+\frac{f_\alpha f_\beta}{\sigma_\Sigma^2}\right| \ d\mu_g\leq c\sigma^{-\frac52-\delta}.
\end{equation}
On the other hand, for $\alpha>3$ we have 
\begin{equation}\label{highmagnDelta}
\lambda_\alpha>\frac{5}{\sigma_\Sigma^2}.
\end{equation}
\end{lem}
\begin{rem} As a byproduct of the proof of \cite[Lemma 2]{cederbaum}, it is important to keep in mind that such an orthonormal system also satisfies the inequality $\|f_i\|_{H^2(\Sigma)}\leq c$, for $i\in\{1,2,3\}$.
\end{rem}
This description of the spectrum of the Laplace-Beltrami operator allows us to define a decomposition of $L^2(\Sigma)$ in terms of the eigenfunctions of $-\Delta$.
\begin{dfn} Let $\Sigma$ be a surface and consider the Hilbert space $L^2(\Sigma)$ equipped with the standard scalar product in $L^2(\Sigma)$. Consider the orthonormal system constructed in Lemma \ref{ortonormalsystemL2}. Then, for every $\textup w\in L^2(\Sigma)$ we define
\begin{equation*}
\textup w^0:=\langle w,f_0\rangle_2f_0=\fint_\Sigma \textup w \ d\mu_g, \qquad \textup w^t:=\sum_{\alpha=1}^3\langle \textup w,f_\alpha\rangle_2f_\alpha.
\end{equation*}
We call $\textup w^0$ the \textit{mean part} of $\textup w$, and $\textup w^t$ the \textit{translational part} of $\textup w$. Finally, we set
\begin{equation*}
\textup w^d:=\textup w-\textup w^t,
\end{equation*}
the so called \textit{difference part}, which obviously also contains the information about $\textup w^0$. Since $\textup w^t$ and $\textup w^d$ are spanned by distinct eigenvalues of $-\Delta$, which are orthonormal, it follows  that $\langle \textup w^t,\textup w^d\rangle_2=0$. When $\textup w^0=0$, we will often call $\textup w^d$ also \textit{deformational part}. 
\end{dfn}
\begin{rem}\begin{enumerate}[label=(\roman*)]
\item First of all, since $\Sigma$ is compact, $f_0\equiv |\Sigma|^{-\frac12}$. It is then clear that if $\textup w$ has zero mean, then $\textup w^0=0$. 
\item The so called \textit{traslational part} of a function $\textup w$ is the projection of $\textup w$ on the finite-dimensional space $\textup{span}\{f_1,f_2,f_3\}$, the eigenfunctions of $-\Delta$ obtained as a perturbation of the coordinate functions on the sphere. 
\end{enumerate}
\end{rem}
\begin{lem}\label{stimaprecisaautolem4}
There exist $c>0$ and $\sigma_0>1$ such that, if $\Sigma\in \overline{\mathcal W}_\sigma^\eta(B_1,B_2)$ and it also is $(\cin,\sigma)$-almost CMC with $\sigma\geq \sigma_0$, the complete orthonormal system in $L^2(\Sigma)$ given by \textup{Lemma \ref{ortonormalsystemL2}} is such that 
\begin{equation}\label{eq490lambdamigl}
\left|\lambda_\alpha-\frac{\hbar^2}{2}-\frac{6m_\mathcal{H}(\Sigma)}{\sigma_\Sigma^3}-\int_\Sigma \overline{\textup{Ric}}(\nu,\nu) f_\alpha^2 \ d\mu_g\right|\leq c\sigma^{-3-\delta}, \qquad \alpha\in\{1,2,3\},
\end{equation}
and the corresponding eigenfunctions $f_1,f_2,f_3$ satisfy
\begin{equation}\label{Ricineq4128}
\left|\int_\Sigma \overline{\textup{Ric}}(\nu,\nu) f_\alpha f_\beta\right|\leq c\sigma^{-3-\delta}, \qquad \alpha\neq \beta, \quad \alpha,\beta\in\{1,2,3\}.
\end{equation}
\end{lem}
\begin{proof} The proof is analogous to the one of \cite[Lemma 3]{cederbaum}. We remark that, using that $|\overset{\circ}{A}|=\textup O(\sigma^{-\frac32-\delta})$, we can write the Gauss equation as 
\begin{equation}\label{gausseq57}
\begin{aligned}
\textup S^\Sigma &=\overline{\textup S}-2\overline{\textup{Ric}}(\nu,\nu)-|\overset{\circ}{A}|^2+\frac{H^2}{2}
\\
&=\overline{\textup{S}}-2\overline{\textup{Ric}}(\nu,\nu)+\frac{\hbar^2}{2}+\frac{(\mathcal H-\hbar)^2}{2}+\hbar(\mathcal H-\hbar)+\textup O(\sigma^{-3-\delta}).
\end{aligned}
\end{equation}
since, using Lemma \ref{propertiesinitialdata} and the roundness, we find $H-\mathcal H=\textup O(\sigma^{-2-2\delta})$ and $\mathcal H=\textup O(\sigma^{-1})$. We thus set $\mathcal R:=\frac{(\mathcal H-\hbar)^2}{2}+\hbar(\mathcal H-\hbar)+\textup O(\sigma^{-3-\delta})$ and, analogously to the proof of \cite[Prop. 3.6]{vpmcf}, we get
\begin{equation}
\begin{aligned}
&\left|\lambda_\alpha^2\delta_{\alpha\beta}-\int_\Sigma \textup S^\Sigma \langle \nabla f_\alpha,\nabla f_\beta\rangle \ d\mu_g\right|
\\
=&\left|\lambda_\alpha^2\delta_{\alpha\beta}-\int_\Sigma \left(\left(\overline{\textup S}-2\overline{\textup{Ric}}(\nu,\nu)\right)+\left(\frac{\hbar^2}{2}+\mathcal R\right)\right) \langle \nabla f_\alpha,\nabla f_\beta\rangle \ d\mu_g\right|.
\end{aligned}
\end{equation}
Since the spacelike case corresponds to $\mathcal R\equiv 0$, in the spacetime case we just have to estimate 
\begin{equation}\label{adproofspt}
\left|\int_\Sigma \mathcal R\langle\nabla f_\alpha,\nabla f_\beta\rangle \ d\mu\right|=\left|\int_\Sigma \left(\frac{(\mathcal H-\hbar)^2}{2}+\hbar(\mathcal H-\hbar)+\textup O(\sigma^{-3-\delta})\right)\langle\nabla f_\alpha,\nabla f_\beta\rangle \ d\mu\right|.
\end{equation}
Remember, comparing this proof with the one of \cite[Prop. 3.6]{vpmcf}, that the aim is to show that this remainder is of order $\textup O(\sigma^{-5-\delta})$. Since $\|f_\alpha\|_{H^2(\Sigma)}=\textup O(1)$, and so $\|\nabla f_\alpha\|_{L^2(\Sigma)}\leq \frac{C}{\sigma}$, notice that $\displaystyle \left|\int_\Sigma \langle \nabla f_\alpha,\nabla f_\beta\rangle \ d\mu\right|\leq C\sigma^{-2}$, which bounds the latter term in \eqref{adproofspt}. The other two terms can be bounded, using Young's inequality, equation \eqref{almostCMCness} and $\mathcal H-H=\textup O(\sigma^{-2-2\delta})$, by
\begin{equation}
C\sigma^{-1}\|\mathcal H-\hbar\|_{L^2}\|\nabla f_\alpha\|_{L^4}\|\nabla f_\beta\|_{L^4}\leq C\sigma^{-5-2\delta}
\end{equation}
using that $\|f_\alpha\|_{W^{1,4}}\leq C\sigma^{-\frac12}\|f_\alpha\|_{H^2}$, and so $\|\nabla f_i\|_4\leq C\sigma^{-\frac32}$. Thus we obtain
\begin{equation}
\left|\left(\lambda_\alpha^2-\frac{\hbar^2}{2}\lambda_\alpha\right)\delta_{\alpha\beta}-\int_\Sigma \left(\overline{\textup S}-2\overline{\textup{Ric}}(\nu,\nu)\right)\left(\frac{3\delta_{\alpha\beta}}{\sigma^2|\Sigma|_g}-\frac{f_\alpha f_\beta}{\sigma^2}\right) \ d\mu_g\right|\leq C\sigma^{-5-\delta},
\end{equation}
which is analogous to inequality obtained in \cite[Prop. 3.6]{vpmcf}.
\end{proof}
\begin{rem}\label{rem2039} Observe that equation \eqref{eq490lambdamigl} gives the following bound. Since $\delta\in (0,\frac12]$, $|m_\mathcal{H}(\Sigma)|\leq 2|E_\textup{ADM}|$, $\overline{\textup{Ric}}_{\vec x}=\textup O(|\vec x|^{-\frac52-\delta})$ and $\|f_\alpha\|_2=1$, then
\begin{equation*}
\left|\lambda_\alpha-\frac{\hbar^2}{2}\right|=\textup O(\sigma^{-\frac52-\delta}),
\end{equation*}
with a constant possibly depending on $|E_\textup{ADM}|$.
\end{rem}
The proof of the following Lemma is similar to \cite[Prop. 2]{cederbaum}.
\begin{lem}\label{behavLem5}
There exist $c>0$ and $\sigma_0>1$ such that, if $\Sigma\in \overline{\mathcal W}_\sigma^\eta(B_1,B_2)$ with $\sigma\geq \sigma_0$, for every $\textup w$, $\textup v\in H^2(\Sigma)$ it holds 
\begin{equation*}
\left|\int_\Sigma (L\textup w^t)\textup v^t \ d\mu-\frac{6m_\mathcal{H}(\Sigma)}{\sigma_\Sigma^3}\int_\Sigma \textup w^t\textup v^t \ d\mu\right|\leq \frac{c\|\textup w\|_{L^2(\Sigma)}\|\textup v\|_{L^2(\Sigma)}}{\sigma^{3+\delta}}.
\end{equation*}
\end{lem}
This leads to the following
\begin{lem}\label{decompstabI}
There exist $c>0$ and $\sigma_0>1$ such that, if $\Sigma$ belongs to $\overline{\mathcal W}_\sigma^\eta(B_1,B_2)$ and it is $(\cin,\sigma)$-almost CMC with $\sigma\geq \sigma_0$, for every $\textup w\in H^2(\Sigma)$ such that $\textup w^0=0$ we find
\begin{enumerate}[label=\textup{(\roman*)}]\item The translational part, in view of \textup{Lemma \ref{behavLem5}}, satisfies
\begin{equation*}
\int_\Sigma (L\textup w^t)\textup w^t \ d\mu\geq \frac{6m_\mathcal H(\Sigma)}{\sigma_\Sigma^3}\|\textup w^t\|_2^2-c\sigma^{-3-\delta}\|\textup w\|_2^2;
\end{equation*}
\item The remaining part satisfies 
\begin{equation*}
\int_\Sigma (L\textup w^d)(\textup w^d) \ d\mu\geq \frac{7}{4\sigma_\Sigma^2}\int_\Sigma (\textup w^d)^2 \ d\mu.
\end{equation*}
\end{enumerate}
\end{lem}
\begin{proof} Point (ii) follows from 
\begin{align*}
\int_\Sigma (L\textup w^d)(\textup w^d) \ d\mu &=\int_\Sigma \textup w^d(-\Delta \textup w^d) \ d\mu\\
&\quad -\int_\Sigma \left(\frac{\hbar^2}{2}+\frac{(\mathcal H-\hbar)^2}{2}+\hbar(\mathcal H-\hbar)+\textup O(\sigma^{-\frac52-\delta})\right)(\textup w^d)^2 \ d\mu.
\end{align*}
Combining this with $\hbar=\frac2{\sigma_\Sigma}+\textup O(\sigma^{-\frac32-\delta})$, together with Lemma \ref{cor1}, i.e. $\|\mathcal H-\hbar\|_{L^\infty(\Sigma)}\leq C\sigma^{-\frac32-\delta}$, and equation \eqref{highmagnDelta}, we get
\begin{equation*}
\int_\Sigma (L\textup w^d)(\textup w^d) \ d\mu\geq \left(\frac{5}{\sigma_\Sigma^2}-\frac{2}{\sigma_\Sigma^2}+\textup O(\sigma^{-\frac52-\delta})\right)\int_\Sigma (\textup w^d)^2 \ d\mu\geq \frac7{4\sigma_\Sigma^2}\int_\Sigma (\textup w^d)^2 \ d\mu,
\end{equation*}
where we also used the equivalence of the radii $\sigma$ and $\sigma_\Sigma$ for surfaces in the class.
\end{proof}
\begin{lem}\label{lemmaintermedio2}
There exist $c>0$ and $\sigma_0>1$ such that, if $\Sigma\in \overline{\mathcal W}_\sigma^\eta(B_1,B_2)$ with $\sigma\geq \sigma_0$, for every $\textup w\in H^2(\Sigma)$ it holds 
\begin{equation*}\|L\textup w^t\|_2^2\leq c\sigma^{-5-2\delta}\|\textup w\|_2^2.
\end{equation*}
\end{lem}
\begin{proof}
We estimate
\begin{equation*}
\|L\textup w^t\|_2\leq \|-\Delta \textup w^t-\frac{\hbar^2}{2}\textup w^t\|_2+\left\|\frac{(\mathcal H-\hbar)^2}{2}\textup w^t+\hbar(\mathcal H-\hbar)\textup w^t\right\|_2+\textup O(\sigma^{-\frac52-\delta})\|\textup w\|_2.
\end{equation*}
Using the definition of $\textup w^t$, and, by Remark \ref{rem2039}, $|\lambda_i-\frac{\hbar^2}{2}|^2=\textup O(\sigma^{-5-2\delta})$, we have
\begin{equation*}
\left\|-\Delta \textup w^t-\frac{\hbar^2}{2}\textup w^t\right\|_2^2\leq C\sigma^{-5-2\delta}\|\textup w\|_2^2.
\end{equation*}
Moreover, we conclude with the estimate
\begin{equation*}
\left\|\frac{(\mathcal H-\hbar)^2}{2}\textup w^t+\hbar(\mathcal H-\hbar)\textup w^t\right\|_2\leq 10\sigma^{-1}\|(\mathcal H-\hbar)\textup w^t\|_2\leq C\sigma^{-\frac52-\delta}\|\textup w\|_2,
\end{equation*}
using again Lemma \ref{cor1}.
\end{proof}
The previous Lemmas lead to the following conclusion.
\begin{prp}\label{spacetimespectrath}
There exist $c>0$ and $\sigma_0>1$ such that, if $\Sigma$ belongs to $\overline{\mathcal W}_\sigma^\eta(B_1,B_2)$ and it is $(\cin,\sigma)$-almost CMC with $\sigma\geq \sigma_0$,
\begin{equation*}
\inf\left\{\int_\Sigma (L\textup w)\textup w \ d\mu: \ \|\textup w\|_{L^2(\Sigma)}=1, \quad \int_\Sigma \textup w \ d\mu=0\right\}\geq \frac{2E_\textup{ADM}}{\sigma_\Sigma^3}.
\end{equation*}
\end{prp}
\begin{proof}
Decomposing the operator $L$ as follows 
\begin{equation*}
\int_\Sigma (L\textup w)\textup w \ d\mu=\int_\Sigma (L\textup w^t)\textup w^t \ d\mu+2\int_\Sigma (L\textup w^t)\textup w^d \ d\mu+\int_\Sigma (L\textup w^d)\textup w^d \ d\mu,
\end{equation*}
and using Lemma \ref{decompstabI}, together with the parametric Young's inequality with $\varepsilon^{-1}=(4\sigma_\Sigma^2)^{-1}$ for the intermediate term, we get 
\begin{equation*}
\begin{aligned}
\int_\Sigma (L\textup w)\textup w \ d\mu\geq & \ \frac{6m_\mathcal H(\Sigma)}{\sigma_\Sigma^3}\|\textup w^t\|_2^2-c\sigma^{-3-\delta}\|\textup w\|_2^2+\frac7{4\sigma_\Sigma^2}\int_\Sigma (\textup w^d)^2 \ d\mu
\\
&-4\sigma_\Sigma^2 \|L\textup w^t\|_2^2-\frac{\|\textup w^d\|_2^2}{4\sigma_\Sigma^2}.
\end{aligned}
\end{equation*}
Using \eqref{mHE43} and choosing $\sigma$ large we have $m_\mathcal{H}(\Sigma)\geq \frac{E_\textup{ADM}}{2}$, and also Lemma \ref{lemmaintermedio2}, we have
\begin{equation*}
\int_\Sigma (L\textup w)\textup w \ d\mu\geq \frac{3E_\textup{ADM}}{\sigma_\Sigma^3}\|\textup w^t\|_2^2-c\sigma^{-3-\delta}\|\textup w\|_2^2+\frac{3}{2\sigma_\Sigma^2}\|\textup w^d\|_2^2.
\end{equation*}
We conclude using $\|\textup w\|_{L^2(\Sigma)}^2=\|\textup w^t\|_{L^2(\Sigma)}^2+\|\textup w^d\|_{L^2(\Sigma)}^2$ choosing $\sigma$ so large that $\displaystyle \frac{3}{2\sigma_\Sigma^2}\geq \frac{3E_\textup{ADM}}{\sigma_\Sigma^3}$.
\end{proof}
\section{Volume preserving spacetime mean curvature flow}\label{vpstmcf4}
\subsection{Definition of the flow and evolution equations}
\begin{dfn} Let $(\M,\overline\g,\overline\K)$ be an initial data set and let $\iota:\Sigma\hookrightarrow \M$ be a closed surface. A time dependent family of immersions $F_t:\Sigma\hookrightarrow\M$, with $t\in [0,T)$ for some $0<T\leq \infty$, which satisfies 
\begin{equation}\label{generalflow513}
\begin{cases}
\frac{\partial}{\partial t} F_t(\cdot)=-\left(\mathcal H(\cdot,t)-\hbar(t)\right)\nu(\cdot,t)\\
F_0=\iota
\end{cases}
\end{equation}
is called a solution to the \textup{volume preserving spacetime mean curvature flow}, with initial value $\iota$.
\end{dfn}
We highlight that the function $\mathcal H$ is an increasing function of the mean curvature $H$. The function $P=g^{ij}\overline\K_{ij}$ in $\mathcal H=\sqrt[q]{H^q-|P|^q}$ depends on the metric induced on $\Sigma_t$, which only involves first order derivatives of the immersion. Thus, without the volume preserving term, the equation is parabolic. However, this term only depends on time, and thus it does not affect the parabolicity and local existence of solutions and uniqueness are ensured.\\
\indent In the following, we will assume that the ambient initial data set is $C_{\frac12+\delta}^2$-asymptotically flat. We write $\Sigma_t:=F_t(\Sigma)$ to denote the solution of the flow at time $t$ and we call $g(t)$ the induced metric and denote by $d\mu_t$ the induced 2-dimensional measure. \\
\indent We recall the evolution equations satisfied by the main geometric quantities on $\Sigma_t$. At each fixed $t$, we choose a frame $\{\vec e_\alpha(t)\}_{\alpha=1}^3$ on $(\M,\overline\g)$ such that  $\{\vec e_1(t),\vec e_2(t)\}$ are tangent vectors on $\Sigma_t$ and $\vec e_3(t):=\nu_t$. The following Lemma collects the equations satisfied by the main geometric quantities on $\Sigma_t$, see \cite{huiskenpolden}.
\begin{lem}\label{evolution41} Let $\{F_t\}_{t\in [0,T)}$ be a solution to the flow \eqref{generalflow513}. Then we have 
\begin{enumerate}[label=\textup{(\roman*)}]
\item $\frac{\partial g_{ij}}{\partial t}=-2\left(\mathcal H-\hbar\right)h_{ij}$;
\item $\frac{\partial}{\partial t}(d\mu_t)=- \left(\mathcal H-\hbar\right)Hd\mu_t$;
\item $\frac{\partial}{\partial t}\nu=\nabla \mathcal H$;
\item $\frac{\partial}{\partial t}h_{ij}=\nabla_i\nabla_j\mathcal H+\left(\mathcal H-\hbar\right)\left(-h_{ik}h_j^k+\overline{\textup{Rm}}_{ikjl}\nu^k\nu^l\right)$;
\item $\frac{\partial H}{\partial t}=\Delta \mathcal H+\left(\mathcal H-\hbar\right)(|A|^2+\overline{\textup{Ric}}(\nu,\nu))$.
\end{enumerate} 
\end{lem}
\noindent \textbf{Notation for the rest of the Section.} In the following, it is convenient to set $\Phi=\Phi(s,\gamma):=\sqrt[q]{s^q-|\gamma|^q}$, so that $\mathcal H=\Phi(H,P)$. We denote by $\Y$ the derivative of $\Phi$ computed with respect to the variable $s$, i.e. $\Y:=\partial_s\Phi\big|_{(s,\gamma)=(H,P)}$. On the other hand, we will denote by $\Psi$ the derivative of $\Phi$ in time, due to the dependence on $P=P(t)$, i.e. $\Psi:=\partial_t\left(\Phi(\rho,P(t))\right)\big|_{\rho=H}$. Thus, 
\begin{equation}\label{decomPsiPhi}
\partial_t\left(\Phi(H,P)\right)=\Y\partial_t H+\Psi.
\end{equation}
This notation is particularly useful since we will mainly take trace of the term $\Y$. In the following, we will have 
\begin{equation}\label{eq561p}
\Upsilon=\left(1-\left(\frac{|P|}{H}\right)^q\right)^{\frac1q-1}, \qquad \nabla\Upsilon=(q-1)\left(\frac{H}{\mathcal H}\right)^{q-2}\left(-\frac{|P|^q}{\mathcal H^2H^{q-1}}\nabla\mathcal H+\frac{1}{\mathcal H}\frac{|P|^{q-2}P}{H^{q-1}}\nabla P\right),
\end{equation}
\begin{equation}\label{estimatePsipre}
|\Psi|=\left|\frac{\partial\Phi(s,P)}{\partial t}\bigg|_{s=H}\right|=\left|\frac{-q\Phi(H,P)}{\left(\Phi(H,P)\right)^q}\frac{(\partial_tP)P}{|P|^{2-q}}\right|\leq C\sigma^{q-1}|P|^{q-1}|\partial_tP|\leq C\sigma^{\frac12-\frac12q-\delta q+\delta}|\partial_tP|,
\end{equation}
where the latter inequality holds assuming \eqref{eq222}, because of Lemma \ref{propertiesinitialdata}. Hypothesis \eqref{eq222} is natural in our setting since we will work solely on round surfaces. Note also that \eqref{eq561p} implies 
\begin{equation}\label{Upsilonineq4545}
|\Upsilon-1|\leq C\left|\frac{P}{H}\right|^q=\textup O(\sigma^{-\frac12q-q\delta}), \qquad \Upsilon-1\geq\underline c\sigma^q|P|^q.
\end{equation} 
\begin{lem} There exists $C>0$ and $\sigma_0>0$ such that, if $\Sigma_t$ satisfies $|A(t)|\leq \sqrt{\frac52}\sigma^{-1}$ and \eqref{eq222} for every $t\in [0,T]$, and $\sigma>\sigma_0$,
\begin{equation}\label{mainpineq59}
|\partial_tP|\leq C\sigma^{-\frac52-\delta}|\mathcal H-\hbar|+C\sigma^{-\frac32-\delta}|\nabla\mathcal H|.
\end{equation}
\end{lem}
\begin{proof} We choose normal coordinates on a point $x^*$ of $\Sigma_{t^*}$, for an arbitrary $t^*\in[0,T]$, say $\{x_1,x_2\}$, and normal coordinates $\{y_1,y_2,y_3\}$ on $y^*:=F_{t^*}(x^*)$ in $\M$. Thus, if $\{\frac{\partial F}{\partial x_i}\}_{i=1}^2$ is the frame induced by the immersion, we notice that 
\begin{equation}
g_{ij}=\left(F^*\overline\g\right)_{ij}=\overline \g\left(\frac{\partial F}{\partial x_i},\frac{\partial F}{\partial x_j}\right), \qquad \overline\K_{ij}=\left(F^*\overline\K\right)_{ij}=\overline\K\left(\frac{\partial F}{\partial x_i},\frac{\partial F}{\partial x_j}\right).
\end{equation}
Thus in particular
\begin{equation}\label{normalcoordin59}
\delta_{\alpha\beta}\frac{\partial F^\alpha}{\partial x_i}\frac{\partial F^\beta}{\partial x_j}=\overline \g\left(\frac{\partial F}{\partial x_i},\frac{\partial F}{\partial x_j}\right)\bigg|_{(x^*,t^*)}=g_{ij}\big|_{(x_*,t_*)}=\delta_{ij}.
\end{equation}
By direct computation, using the symmetry of $\overline \K$ we find that 
\begin{equation}\label{eqevP}
\begin{aligned}
\partial_tP=\partial_t\left(g^{ij}\overline \K_{ij}\right)&=2(\mathcal H-\hbar)h^{ij}\overline \K_{ij}+\overline\nabla_\gamma\overline\K_{ij}\frac{\partial F^\gamma}{\partial t}+2g^{ij}\overline \K\left(\frac{\partial}{\partial t}\left(\frac{\partial F}{\partial x_i}\right),\frac{\partial F}{\partial x_j}\right)\\
&=2(\mathcal H-\hbar)h^{ij}\overline \K_{ij}+\overline\nabla_\gamma\overline\K_{ij}(\hbar-\mathcal H)\nu^\gamma+2g^{ij}\overline \K\left(\frac{\partial}{\partial t}\left(\frac{\partial F}{\partial x_i}\right),\frac{\partial F}{\partial x_j}\right).
\end{aligned}
\end{equation}
Since 
\begin{equation}
\frac{\partial}{\partial t}\left(\frac{\partial F}{\partial x_i}\right)=\frac{\partial}{\partial x_i} \left(\frac{\partial F}{\partial t}\right)=\frac{\partial}{\partial x_i} \left((\hbar-\mathcal H)\nu\right)=-\frac{\partial \mathcal H}{\partial x_i}\nu+(\hbar-\mathcal H)\frac{\partial \nu}{\partial x_i},
\end{equation}
we rewrite \eqref{eqevP} as
\begin{equation}\label{evPestesa512}
\partial_tP=2(\mathcal H-\hbar)h^{ij}\overline \K_{ij}+\overline\nabla_\gamma\overline\K_{ij}(\hbar-\mathcal H)\nu^\gamma-2g^{ij}\overline\K_{\alpha\beta}\nu^\alpha\frac{\partial\mathcal H}{\partial x_i}\frac{\partial F^\beta}{\partial x_j}+2(\hbar-\mathcal H)g^{ij}\overline\K_{\alpha\beta}\frac{\partial \nu^\alpha}{\partial x_i}\frac{\partial F^\beta}{\partial x_j}.
\end{equation}
Note that, in normal coordinates, the Weingarten equation takes the form 
\begin{equation}\label{weingarten513}
\frac{\partial \nu^\alpha}{\partial x_i}\bigg|_{(x^*,t^*)}=h_i^j(x^*,t^*)\frac{\partial F^\alpha}{\partial x_i}\bigg|_{(x^*,t^*)},
\end{equation}
see \cite[Pg. 63]{huiskenpolden}. Thus, computing \eqref{evPestesa512} in the point $(x^*,t^*)$, and estimating, we get
\begin{equation}
|\partial_tP|\leq C|\mathcal H-\hbar||A||\overline\K|+|\overline\nabla \ \overline\K||\mathcal H-\hbar|+C|\overline\K||\nabla \mathcal H|,
\end{equation}
where we used \eqref{weingarten513} combined with \eqref{normalcoordin59} in order to estimate the latter term in \eqref{evPestesa512}. We conclude using that, thanks to the assumption on $|A(t)|$ and \eqref{eq222}, Lemma \ref{propertiesinitialdata} implies that $|\overline\K|\leq C\sigma^{-\frac32-\delta}$ and $|A||\overline\K|+|\overline\nabla \ \overline\K|\leq C\sigma^{-\frac52-\delta}$.
\end{proof}
The $\Phi$-notation, together with helping us avoiding huge formulas in the following, highlights that existence and convergence of the flow could be studied in the case of more general speed functions. However, we just focus our attention on the spacetime flow. We also define $\alpha:(0,1)\to\R$ to be $\alpha(\rho):=\sqrt[q]{1-\rho^q}$, so that $\Phi(s,\gamma)=s\alpha\left(\frac{|\gamma|}{s}\right)$.
\begin{lem}\label{stevolequations52} Along a solution of the volume preserving spacetime mean curvature flow we have 
\begin{equation}\label{eq84v}
\begin{aligned}
\frac{\partial}{\partial t}|\overset{\circ}{A}|^2 = & \ \Delta |\overset{\circ}{A}|^2-2|\nabla\overset{\circ}{A}|^2+\frac{2\hbar}{H}\left\{|A|^4-H\textup{tr}(A^3)\right\}+2|A|^2\left(\frac{H-\hbar}{H}\right)|\overset{\circ}{A}|^2
\\
&+2(\mathcal H-\hbar)\overset{\circ}{h}_{ij}\overline{\textup{Rm}}_{kilj}\nu^k\nu^l-2\left(h_i^l\overline{\textup{Rm}}_{kjkl}+h^{lk}\overline{\textup{Rm}}_{lijk}\right) h_{ij}
\\
&-2\left(\nabla_j\left(\overline{\textup{Ric}}_{i\varepsilon}\nu^\varepsilon\right)+\nabla_l \left(\overline{\textup{Rm}}_{\varepsilon ijl}\nu^\varepsilon\right)\right)\overset{\circ}{h}_{ij}+2|A|^2\left(\frac{H^2-H\mathcal H}{2}\right)
\\
&+2(\mathcal H-H)\textup{tr}(A^3)+\langle \mathcal T,\overset{\circ}{A}\rangle;
\end{aligned}
\end{equation}
\begin{equation}\label{eq101nbis}
\begin{aligned}
\frac{\partial}{\partial t}|\nabla \mathcal H|^2 = & \ \Delta|\nabla \mathcal H|^2-2|\nabla^2\mathcal H|^2+2(\mathcal H-\hbar)h^{ij}\nabla_i\mathcal H\nabla_j\mathcal H+2(|A|^2+\overline{\textup{Ric}}(\nu,\nu))|\nabla \mathcal H|^2
\\
&-2\textup{Ric}^\Sigma(\nabla \mathcal H,\nabla \mathcal H)+2(\mathcal H-\hbar)\langle\nabla|A|^2,\nabla \mathcal H\rangle+2(\mathcal H-\hbar)\langle\nabla\left(\overline{\textup{Ric}}(\nu,\nu)\right),\nabla \mathcal H\rangle
\\
&+2g^{ij}\nabla_i\left(\left(\Y-1\right)\left(\Delta\mathcal H+(\mathcal H-\hbar)(|A|^2+\overline{\textup{Ric}}(\nu,\nu))\right)\right)\nabla_j\mathcal H+2g^{ij}\nabla_i\mathcal H\nabla_j \mathcal H
\end{aligned}
\end{equation}
where $\textup{Ric}^\Sigma$ is the Ricci tensor of $\Sigma$, $\beta:=\alpha-1$ and $\mathcal T:=\left(T_{ij}\right)$ is the tensor defined by
\begin{equation}\label{Temindertens}
\begin{aligned}
T_{ij}:=&\left(\nabla_i\nabla_j H\right)\beta\left(\frac{|P|}{H}\right)+\nabla_j H\beta'\left(\frac{|P|}{H}\right)\nabla_i\left(\frac{|P|}{H}\right)+\nabla_i H\beta'\left(\frac{|P|}{H}\right)\nabla_j\left(\frac{|P|}{H}\right)
\\
&+H\beta''\left(\frac{|P|}{H}\right)\nabla_i\left(\frac{|P|}{H}\right)\nabla_j\left(\frac{|P|}{H}\right)+H\beta'\left(\frac{|P|}{H}\right)\nabla_i\nabla_j\left(\frac{|P|}{H}\right).
\end{aligned}
\end{equation}
\end{lem}
The proof is standard, and it mainly relies on the computations in \cite{huisken1987} and \cite{huiskenpolden}. See moreover \cite{vpmcf}. Observe that the tensor $T$ is the remainder of the Hessian of the function $\Phi$, which, due to the introduction of the auxiliary functions $\alpha$ and $\beta$, is given by $\textup{Hess}(H)$ plus the tensor $T$. Finally, an easy computation shows that, since $q\geq2$,
\begin{equation}\label{controlbeta2}
|\beta(\rho)|\leq c_q\rho^2, \qquad |\beta'(\rho)|\leq c_q\rho, \qquad |\beta''(\rho)|\leq c_q,
\end{equation}
for $\rho<<1$, which is the case we are interested in, since $\rho\sim \frac{|P|}{H}$ which is small on a round surface.
\begin{proof} Using Lemma \ref{evolution41}, we get
\begin{equation}\label{eq525}
\frac{\partial}{\partial t}h_{ij}=\nabla_i\nabla_j\left(\Phi(H,P)\right)+\left(\Phi(H,P)-\hbar\right)\left(-h_{ik}h_j^k+\overline{\textup{Riem}}_{i3j3}\right).
\end{equation}
By $\Phi(H,P)=H\alpha\left(\frac{|P|}{H}\right)$, we have
\begin{equation}
\begin{aligned}
\nabla_i\nabla_j\left(\Phi(H,\cdot)\right)= &\ \left(\nabla_i\nabla_jH\right)\alpha\left(\frac{|P|}{H}\right)\\
&\ +\nabla_j H\alpha'\left(\frac{|P|}{H}\right)\nabla_i\left(\frac{|P|}{H}\right)\\
&\ +\nabla_i H\alpha'\left(\frac{|P|}{H}\right)\nabla_j\left(\frac{|P|}{H}\right)\\
&\ +H\alpha''\left(\frac{|P|}{H}\right)\nabla_i\left(\frac{|P|}{H}\right)\nabla_j\left(\frac{|P|}{H}\right)\\
&\ +H\alpha'\left(\frac{|P|}{H}\right)\nabla_i\nabla_j\left(\frac{|P|}{H}\right)
\end{aligned}
\end{equation}
We moreover define $\beta$ as above, obtaining $\beta'=\alpha'$ and $\beta''=\alpha''$. We thus get
\begin{equation}\label{eq18}
\nabla_i\nabla_j\left(\Phi(H,P)\right)=\nabla_i\nabla_j H+T_{ij}.
\end{equation}
Then \eqref{eq525} becomes
\begin{equation}
\begin{aligned}
\frac{\partial}{\partial t}h_{ij}=& \ \nabla_i\nabla_jH+\left(\Phi(H,\cdot)-\hbar\right)\left(-h_{ik}h_j^k+\overline{\textup{Riem}}_{i3j3}\right)+T_{ij}
\\
=& \ \Delta h_{ij}-Hh_i^lh_{lj}+|A|^2h_{ij}+\left(\Phi(H,\cdot)-\hbar\right)\left(-h_{ik}h_j^k+\overline{\textup{Riem}}_{i3j3}\right)
\\
&-h_i^l\overline{\textup{Rm}}_{kjkl}-h^{lk}\overline{\textup{Rm}}_{lijk}-\nabla_j\left(\overline{\textup{Ric}}_{i\varepsilon}\nu^\varepsilon\right)-\nabla^l\left(\overline{\textup{Rm}}_{\varepsilon ijl}\nu^\varepsilon\right)+T_{ij}
\end{aligned}
\end{equation}
The conclusion follows remarking that 
\begin{equation}
\begin{aligned}
&\Delta|\overset{\circ}{A}|^2-2|\nabla \overset{\circ}{A}|^2+2|A|^2\left(|A|^2-\frac{H\Phi(H,P)}2\right)-2\hbar\textup{tr}(A^3)+H|A|^2\hbar\\
&+2(\Phi(H,P)-H)\textup{tr}(A^3)+\frac{2\hbar}{H}|A|^4-\frac{2\hbar}{H}|A|^4
\\
= \ & \Delta |\overset{\circ}{A}|^2-2|\nabla\overset{\circ}{A}|^2+\frac{2\hbar}{H}|A|^4-2\hbar\textup{tr}(A^3)+2|A|^2\left(1-\frac{\hbar}{H}\right)\left(|A|^2-\frac{H^2}2\right)
\\
&+2|A|^2\left(\frac{H^2}{2}-\frac{H\Phi(H,P)}{2}\right)+2(\Phi(H,P)-H)\textup{tr}(A^3).
\end{aligned}
\end{equation}
Finally, equation \eqref{eq101nbis} follows from Lemma \ref{evolution41} and the Bochner formula.
\end{proof}
\subsection{Evolution of integral quantities}\label{section42}
We now study the evolution of some integral quantities along the flow. Throughout the subsection, $F_t:\Sigma\hookrightarrow \M$ will be a solution to the volume preserving spacetime mean curvature flow \eqref{generalflow513}, in a constrained initial data set $(\M,\overline\g,\overline\K)$, with $t\in [0,T]$ for some $T>0$. We will assume that the surfaces $\Sigma_t$ satisfy properties \eqref{radius1} and \eqref{radius2} of round surfaces for some given suitably large radius $\sigma$. In some results, we further assume 
\begin{equation}\label{Binftycinftyst}
\|H-h\|_{L^\infty(\Sigma_t)}\leq B_\infty \sigma^{-\frac32-\delta}, \quad \left\|\overset{\circ}{A}(t)\right\|_{L^\infty(\Sigma)}\leq B_\infty\sigma^{-\frac32-\delta},
\end{equation}
which are properties satisfied by round surfaces, see Lemma \ref{cor1}, and also
\begin{equation}\label{H^1est411}
\|\mathcal H-\hbar\|_{H^1(\Sigma_t)}\leq \cin\sigma^{-\frac{q}{2}-q\delta}.
\end{equation}
We do not assume apriori that $\Sigma_t$ satisfy properties \eqref{cond2defroundst} and \eqref{cond3defroundst}. We want to analyze the invariance of these properties along the flow. We start estimating the $L^4$ norm of the traceless second fundamental form of $\Sigma_t$. In this result, hypothesis \eqref{Binftycinftyst} is replaced by the milder assumption \eqref{milderassumotion}.
\begin{prp}\label{eqevApallST} Let $\{F_t\}_{t\in [0,T]}$ be a solution to the flow satisfying \eqref{radius1} and \eqref{radius2}. Suppose in addition
\begin{equation}\label{milderassumotion}
\|H-h\|_{L^\infty(\Sigma_t)}\leq \frac{1}{20\sigma};
\end{equation}
Then there exist a constant $C=C(\overline c,\delta)>0$ and a radius $\sigma_0=\sigma_0(\delta,\overline c)>0$ such that if $\sigma>\sigma_0$ then 
\begin{equation}\label{evL4_52st}
\frac{\textup{d}}{\textup{d}t}\int_\Sigma |\overset{\circ}{A}|^4 \ d\mu_t\leq-2\int_\Sigma |\overset{\circ}{A}|^{2}|\nabla \overset{\circ}{A}|^2 \ d\mu_t
 -\frac{1}{2\sigma^2}\int_\Sigma |\overset{\circ}{A}|^4 \ d\mu_t+C\sigma^{-6-4\delta}.
\end{equation} 
As a consequence, if $\displaystyle\int_\Sigma |\overset{\circ}{A}|^4 \ d\mu_0<B_1\sigma^{-4-4\delta}$ and $B_1>2C$, then $\displaystyle\int_\Sigma |\overset{\circ}{A}|^4 \ d\mu_t<B_1\sigma^{-4-4\delta}$ for every $t\in [0,T]$.
\end{prp}
\begin{proof} The proof consists in studying the evolution of the quantity $\|\overset{\circ}{A}\|_{L^4(\Sigma_t)}^4$. Notice that, using Lemma \ref{stevolequations52} and comparing the evolution \eqref{eq84v} with the one in \cite[Prop. 3.6]{vpmcf}, we only have to study the "new" term
\begin{equation*}
\int_\Sigma \langle \mathcal T,\overset{\circ}{A}\rangle|\overset{\circ}{A}|^{2} \ d\mu.
\end{equation*}
Once this term have been estimated, we can conclude the proof following the one of \cite[Prop. 3.6]{vpmcf} combined with Lemma \ref{propertiesinitialdata} in order to deal with the difference between $H$ and $\mathcal H$. Thus, it only remains to prove that there exists $C$ such that, for every fixed $\varepsilon>0$ and $\sigma_0$ suitably large it holds
\begin{equation}\label{proofappendix}
\left|\int_\Sigma \langle \mathcal T,\overset{\circ}{A}\rangle|\overset{\circ}{A}|^{2} \ d\mu\right|\leq \frac{\varepsilon}{\sigma^2}\int_\Sigma |\overset{\circ}{A}|^4 \ d\mu+\varepsilon\int_\Sigma |\overset{\circ}{A}|^{2}|\nabla\overset{\circ}{A}|^2 \ d\mu+C\sigma^{-6-4\delta}.
\end{equation}
Multiplying equation \eqref{Temindertens} by $\overset{\circ}{h}_{ij}|\overset{\circ}{A}|^2$ and integrating we get 
\begin{equation*}
\begin{aligned}
\int_\Sigma T_{ij}\overset{\circ}{h}_{ij}|\overset{\circ}{A}|^2 \ d\mu
&= \int_\Sigma \Big\{\nabla_i\nabla_j H\beta\left(\frac{|P|}{H}\right)+\nabla_j H\beta'\left(\frac{|P|}{H}\right)\nabla_i\left(\frac{|P|}{H}\right)\\
&\quad +\nabla_i H\beta'\left(\frac{|P|}{H}\right)\nabla_j\left(\frac{|P|}{H}\right)+H\beta''\left(\frac{|P|}{H}\right)\nabla_i\left(\frac{|P|}{H}\right)\nabla_j\left(\frac{|P|}{H}\right)\\
&\quad +H\beta'\left(\frac{|P|}{H}\right)\nabla_i\nabla_j\left(\frac{|P|}{H}\right)\Big\}\overset{\circ}{h}_{ij} |\overset{\circ}{A}|^2 \ d\mu.
\end{aligned}
\end{equation*}
Integration by parts, Lemma \ref{propertiesinitialdata} and $\nabla_i\left(\overset{\circ}{h_{ij}}|\overset{\circ}{A}|^{2}\right)=\nabla_i\overset{\circ}{h_{ij}}|\overset{\circ}{A}|^{2}+2\overset{\circ}{h_{ij}}|\overset{\circ}{A}|\nabla_i|\overset{\circ}{A}|$ imply 
\begin{equation}\label{eq328}
\begin{aligned}
\left|\int_\Sigma (\nabla_i\nabla_j H)\beta\left(\frac{P}{H}\right)\overset{\circ}{h_{ij}}|\overset{\circ}{A}|^{2}\ d\mu\right|&\leq C\sigma^{-\frac12-\delta}\int_\Sigma \nabla_j H\left|\frac{(\nabla_i P)H-P\nabla_i H}{H^2}\right||\overset{\circ}{A}|^{3} \ d\mu
\\
&\quad +C\sigma^{-1-2\delta}\int_\Sigma |\nabla H||\overset{\circ}{A}|^{2}|\nabla \overset{\circ}{A}| \ d\mu.
\end{aligned}
\end{equation}
Using again Lemma \ref{propertiesinitialdata}, the parametric Young's inequality, $|\nabla H|^2\leq C|\nabla\overset{\circ}{A}|^2+C|\overline{\textup{Ric}}|^2$ (see \cite{huisken86}) and \eqref{asymptc8} we get 
\begin{equation}\label{eqA4}
\left|\int_\Sigma (\nabla_i\nabla_j H)\beta\left(\frac{P}{H}\right)\overset{\circ}{h_{ij}}|\overset{\circ}{A}|^{2}\ d\mu\right|\leq \frac{\varepsilon}{\sigma^2}\int_\Sigma |\overset{\circ}{A}|^4 \ d\mu+\varepsilon\int_\Sigma |\overset{\circ}{A}|^{2}|\nabla\overset{\circ}{A}|^2 \ d\mu+C\sigma^{-6-4\delta}.
\end{equation}
Moreover, we estimate 
\begin{equation*}
\begin{aligned}
\left|\int_\Sigma (\nabla_i H)\beta'\left(\frac{P}{H}\right)\nabla_j\left(\frac{P}{H}\right)\overset{\circ}{h}_{ij}|\overset{\circ}{A}|^{2} \ d\mu\right|\leq & \  C\sigma^{-2-2\delta}\int_\Sigma \left(|\nabla H||\overset{\circ}{A}|\right)|\overset{\circ}{A}|^{2} \ d\mu\\
&+C\sigma^{-2\delta}\int_\Sigma |\nabla H|^2|\overset{\circ}{A}|^{3} \ d\mu.
\end{aligned}
\end{equation*}
The second addend can be estimated as in \eqref{eqA4}, while the first addend, using Young's inequality, is bounded by
\begin{equation}\label{auxineq6109}
C\sigma^{-2-2\delta}\int_\Sigma \left(\frac{|\nabla H|^2}{2}+\frac{|\overset{\circ}{A}|^2}{2}\right)|\overset{\circ}{A}|^{2} \ d\mu.
\end{equation}
Again, the first addend of \eqref{auxineq6109} can be treated as in \eqref{eqA4}, for $\sigma$ large.\\
\indent We can also bound the term
\begin{equation*}
\begin{aligned}
\left|\int_\Sigma H\beta''\left(\frac{P}{H}\right)\nabla_i\left(\frac{P}{H}\right)\nabla_j\left(\frac{P}{H}\right)\overset{\circ}{h_{ij}}|\overset{\circ}{A}|^{2} \ d\mu\right|\leq & \ C\sigma^{-2\delta}\int_\Sigma |\nabla H|^2|\overset{\circ}{A}|^{2}|\overset{\circ}{A}|\ d\mu\\
&+ C\sigma^{-4-2\delta}\int_\Sigma |\overset{\circ}{A}|^{3} \ d\mu
\\
&+2\int_\Sigma \frac{|P||\nabla H||\nabla P|}{H^2}|\overset{\circ}{A}|^{3} \ d\mu.
\end{aligned}
\end{equation*}
We conclude as in \eqref{eqA4}, also using $|\overset{\circ}{A}|\leq C\sigma^{-\frac12}$, for $\sigma$ large and Young's inequality.\\
\indent Finally, integrating by parts and using the decay of $\beta$ we get
\begin{equation}\label{addend63}
\begin{aligned}
&\left|\int_\Sigma H \, \beta'\left(\frac{P}{H}\right)\nabla_i\nabla_j\left(\frac{P}{H}\right) \overset{\circ}{h_{ij}} |\overset{\circ}{A}|^{2} \, d\mu\right|\\
&\qquad\qquad\qquad\leq C\sigma^{-\frac12-\delta}\int_\Sigma |\nabla H|\, \left|\frac{(\nabla P)H-(\nabla H)P}{H^2}\right||\overset{\circ}{A}|^{3} \, d\mu
\\
&\qquad\qquad\qquad\quad +C\int_\Sigma H \, \left|\frac{(\nabla P)H-(\nabla H)P}{H^2}\right|^2|\overset{\circ}{A}|^{3} \, d\mu
\\
&\qquad\qquad\qquad\quad +C\sigma^{-\frac12-\delta}\int_\Sigma H \, \left|\frac{(\nabla P)H-(\nabla H)P}{H^2}\right||\nabla \overset{\circ}{A}| \, |\overset{\circ}{A}|^{2} \, d\mu.
\end{aligned}
\end{equation}
The first addend can be dealt with as \eqref{eqA4}, while the second as in \eqref{auxineq6109}. The third addend in \eqref{addend63} is bounded by
\begin{equation*}
C\sigma^{-\frac12-\delta}\int_\Sigma |\nabla P|\left|\nabla |\overset{\circ}{A}|\right||\overset{\circ}{A}|^{2} \ d\mu+C\sigma^{-1-2\delta}\int_\Sigma |\nabla H|\left|\nabla |\overset{\circ}{A}|\right||\overset{\circ}{A}|^{2} \ d\mu,
\end{equation*}
and we conclude by Young's inequality, combined with the estimate $|\nabla H|^2\leq C|\nabla\overset{\circ}{A}|^2+C|\overline{\textup{Ric}}|^2$.
\end{proof}
We next estimate the rate of change of the volume preserving term $\hbar(t)$ and of the $L^4$ norm of $\mathcal H-\hbar$. In particular the following Lemma employs some techniques learned in \cite{li}.
\begin{lem}\label{rateh(t)st} Let $(\Sigma,F_t)$ be a solution to the volume preserving spacetime mean curvature flow for $t\in [0,T]$, satisfying properties \eqref{radius1}, \eqref{radius2}, \eqref{Binftycinftyst} and \eqref{H^1est411}. Then, there exist $C=C(\overline c)>0$ and $\sigma_0=\sigma_0(\overline c)>1$, such that, if $\sigma>\sigma_0$,
\begin{equation}\label{roughestL2t}
\frac{\textup d}{\textup dt}\int_\Sigma (\mathcal H-\hbar)^2 \ d\mu_t\leq -\frac12\int_\Sigma |\nabla\mathcal H|^2 \ d\mu_t+C\sigma^{-2}\int_\Sigma (\mathcal H-\hbar)^2 \ d\mu_t.
\end{equation}
Moreover, there exists a constant $c=c(\cin,\overline c)>0$ and a universal constant $C=C(\overline c)>0$ such that 
\begin{equation}\label{rateeqref}
|\dot \hbar(t)|\leq c \sigma^{-\frac72-\frac12q-\delta-\delta q},
\end{equation}
\begin{equation*}\label{eq41645}
\frac{\textup{d}}{\textup{d}t}\int_\Sigma (\mathcal H-\hbar)^4 \ d\mu_t\leq -12\int_\Sigma (\mathcal H-\hbar)^2|\nabla \mathcal H|^2 \ d\mu_t+C\sigma^{-2}\int_\Sigma (\mathcal H-\hbar)^4 \ d\mu_t+cB_\infty \sigma^{-5-\frac32 q-2\delta-3\delta q},
\end{equation*}
provided $\sigma\geq \sigma_0$, for a suitably $\sigma_0=\sigma_0(B_\infty,\cin,\overline c,\delta)$.
\end{lem}
\begin{proof} We first prove \eqref{rateeqref}. By definition of $\hbar$, we get
\begin{equation*}
\begin{aligned}
&|\Sigma_t|\dot \hbar(t) \\
&\qquad =\int_\Sigma \left(\frac{\partial H}{\partial t}+\left(\Upsilon-1\right)\frac{\partial H}{\partial t}+\Psi\right) \ d\mu_t+\int_\Sigma \mathcal HH(\hbar-\mathcal H) \ d\mu_t-\hbar\int_\Sigma H(\hbar-\mathcal H) \ d\mu_t\\
&\qquad =\int_\Sigma (\mathcal H-\hbar)\left(|\overset{\circ}{A}|^2+\overline{\textup{Ric}}(\nu,\nu)\right) \ d\mu_t+\int_\Sigma \left(\frac{H^2}{2}-\mathcal H H\right)\left(\mathcal H-\hbar\right) \ d\mu_t\\
&\qquad  \quad-\hbar\int_\Sigma H(\hbar-\mathcal H) \ d\mu_t+\int_\Sigma \left((\Upsilon-1)\frac{\partial H}{\partial t}+\Psi\right) \ d\mu_t\\
&\qquad =\int_\Sigma (\mathcal H-\hbar)\left(|\overset{\circ}{A}|^2+\overline{\textup{Ric}}(\nu,\nu)\right) \ d\mu_t-\int_\Sigma \frac{H^2}{2}\left(\mathcal H-\hbar\right) \ d\mu_t\\
&\qquad \quad-\int_\Sigma (\mathcal H-H)H(\mathcal H-\hbar) \ d\mu_t-\hbar\int_\Sigma H(\hbar-\mathcal H) \ d\mu_t+\int_\Sigma \left((\Upsilon-1)\frac{\partial H}{\partial t}+\Psi\right) \ d\mu_t\\
&\qquad=\int_\Sigma (\mathcal H-\hbar)\left(|\overset{\circ}{A}|^2+\overline{\textup{Ric}}(\nu,\nu)\right) \ d\mu_t-\int_\Sigma \frac{\mathcal H^2}{2}\left(\mathcal H-\hbar\right) \ d\mu_t \\
&\qquad \quad +\int_\Sigma \left(\frac{\mathcal H^2-H^2}{2}\right)(\mathcal H-\hbar) \ d\mu_t-\int_\Sigma (\mathcal H-H)H(\mathcal H-\hbar) \ d\mu_t-\hbar\int_\Sigma H(\hbar-\mathcal H) \ d\mu_t\\
&\qquad \quad +\int_\Sigma \left((\Upsilon-1)\frac{\partial H}{\partial t}+\Psi\right) \ d\mu_t
\end{aligned}
\end{equation*}
To estimate the above terms, we first note that 
\begin{equation}
\left|\int_\Sigma (\mathcal H-\hbar)\left(|\overset{\circ}{A}|^2+\overline{\textup{Ric}}(\nu,\nu)\right) \ d\mu_t\right|\leq c\sigma^{-\frac32-\frac{q}2-q\delta},
\end{equation}
and that, in addition, the following identity holds
\begin{equation}
\begin{aligned}
&-\int_\Sigma \frac{\mathcal H^2}{2}\left(\mathcal H-\hbar\right) \ d\mu_t-\hbar\int_\Sigma H(\hbar-\mathcal H) \ d\mu_t\\
= &-\int_\Sigma \frac{\mathcal H^2}{2}\left(\mathcal H-\hbar\right) \ d\mu_t+\hbar\int_\Sigma \left(\mathcal H-\hbar\right)^2 \ d\mu_t-\hbar\int_\Sigma \left\{(H-h)-(\mathcal H-\hbar)\right\}\left(\hbar-\mathcal H\right) \ d\mu_t\\
= & -\frac12\int_\Sigma(\mathcal H-\hbar)^3 \ d\mu_t-\hbar\int_\Sigma \left\{(H-h)-(\mathcal H-\hbar)\right\}\left(\hbar-\mathcal H\right) \ d\mu_t,
\end{aligned}
\end{equation}
using also $\displaystyle\int_\Sigma (\mathcal H-\hbar) \ d\mu_t=0$, where, thanks to Lemma \ref{propertiesinitialdata},
\begin{equation}\label{eq419}
\left|\hbar\int_\Sigma \left\{(H-h)-(\mathcal H-\hbar)\right\}\left(\hbar-\mathcal H\right) \ d\mu_t\right|\leq c\sigma^{-1-q-2q\delta}.
\end{equation}
Since the remaining addend can be estimated in a similar way to \eqref{eq419}, we get 
\begin{equation}
|\Sigma_t||\dot\hbar(t)|\leq \left|\frac12\int_\Sigma (\mathcal H-\hbar)^3 \ d\mu_t\right|+c\sigma^{-\frac32-\frac{q}2-q\delta}+\left|\int_\Sigma\left((\Upsilon-1)\frac{\partial H}{\partial t}+\Psi\right) \ d\mu_t\right|
\end{equation}
Observe that the term $\displaystyle \int_\Sigma (\mathcal H-\hbar)^3 \ d\mu_t$ can be easily bounded using \eqref{Binftycinftyst} and \eqref{H^1est411}. Finally, we estimate
\begin{align*}
\int_\Sigma \left((\Upsilon-1)\frac{\partial H}{\partial t}+\Psi\right) \ d\mu_t & = - \int_\Sigma \nabla \Upsilon\cdot\nabla\mathcal H \ d\mu_t\\
&\quad +\int_\Sigma (\Upsilon-1)\left(|A|^2+\overline{\textup{Ric}}(\nu,\nu)\right)(\mathcal H-\hbar) \ d\mu_t\\
&\quad +\int_\Sigma \Psi \ d\mu_t.
\end{align*}
To estimate this term, we observe that equation \eqref{eq561p}, together with the inequalities
\begin{equation*}
\left|\frac{P^p}{\mathcal H^2H^{p-1}}\right|\leq c_q\sigma^{1-\frac12q-q\delta}, \qquad \left|\frac{1}{\mathcal H}\frac{|P|^{q-2}P}{H^{q-1}}\right|\leq c_q\sigma^{\frac32-\frac12q-\delta(q-1)},
\end{equation*}
imply
\begin{equation*}
\left|\int_\Sigma \nabla\Upsilon \cdot \nabla\mathcal H \ d\mu_t\right|\leq c\sigma^{-\frac32-\frac{q}2-\delta q-\delta}.
\end{equation*}
Similarly, we also have 
\begin{equation}
\left|\int_\Sigma (\Upsilon-1)\left(|A|^2+\overline{\textup{Ric}}(\nu,\nu)\right)(\mathcal H-\hbar) \ d\mu_t\right|\leq c\sigma^{-1-q-q\delta}.
\end{equation}
We conclude combining \eqref{estimatePsipre}, \eqref{mainpineq59} and assumption \eqref{H^1est411}, in order to get
\begin{equation}\label{Psi4equalc}
\int_\Sigma \left|\Psi\right| \ d\mu_t \leq c\sigma^{-1-\frac12q-\delta q}\|\mathcal H-\hbar\|_{L^2(\Sigma_t)}+c\sigma^{-\frac12q-\delta q}\|\nabla\mathcal H\|_{L^2(\Sigma_t)}\leq c\sigma^{-1-q-2\delta q}.
\end{equation}
Equation \eqref{rateeqref} follows dividing by $|\Sigma_t|\geq (7/2)\pi\sigma^2$.\\
\indent We now prove \eqref{roughestL2t}. We compute the evolution of $\|\mathcal H-\hbar\|_{L^2(\Sigma,\mu_t)}^2$, obtaining
\begin{equation*}
\begin{aligned}
\frac{\textup d}{\textup{d}t}\int_\Sigma (\mathcal H-\hbar)^2 \ d\mu_t &=2\int_\Sigma (\mathcal H-\hbar)\left(\Delta\mathcal H+(\mathcal H-\hbar)(|A|^2+\overline{\textup{Ric}}(\nu,\nu))\right) \, d\mu_t\\
&\quad + 2\int_\Sigma (\mathcal H-\hbar)(\Upsilon-1)\left(\Delta\mathcal H +(\mathcal H-\hbar)(|A|^2+\overline{\textup{Ric}}(\nu,\nu))\right) \, d\mu_t\\
&\quad + 2\int_\Sigma \Psi(\mathcal H-\hbar) \ d\mu_t-\int_\Sigma H(\mathcal H-\hbar)^3 \, d\mu_t.
\end{aligned}
\end{equation*}
Using integration by parts, the estimate $H+|H-h|+|A|\leq C\sigma^{-1}$, and the inequality
\begin{equation}\label{ineq122}
\begin{aligned}
\left|\int_\Sigma (\mathcal H-\hbar)(\Upsilon-1)\Delta \mathcal H \ d\mu_t\right|=&\left|-\int_\Sigma (\Upsilon-1)|\nabla\mathcal H|^2 \ d\mu_t-\int_\Sigma (\mathcal H-\hbar)\nabla\Upsilon\nabla\mathcal H \ d\mu_t\right|
\\
\leq & \ \varepsilon\int_\Sigma |\nabla\mathcal H|^2 \ d\mu_t+C\int_\Sigma \sigma^{-\frac12q+1-q\delta}|\mathcal H-\hbar||\nabla\mathcal H|^2 \ d\mu_t
\\
&+C\int_\Sigma \sigma^{-\frac12q-1-q\delta}|\mathcal H-\hbar||\nabla\mathcal H| \ d\mu_t
\end{aligned}
\end{equation}
and Young's inequality we get
\begin{equation}
\displaystyle \frac{\textup d}{\textup{d}t}\int_\Sigma (\mathcal H-\hbar)^2 \ d\mu_t\leq-(2-2\varepsilon)\int_\Sigma |\nabla \mathcal H|^2 \ d\mu_t +C\sigma^{-2}\int_\Sigma (\mathcal H-\hbar)^2 \ d\mu_t,
\end{equation}
where we estimated $\displaystyle \int_\Sigma \Psi(\mathcal H-\hbar) \ d\mu_t$ combining \eqref{estimatePsipre} and \eqref{mainpineq59}, i.e. 
\begin{equation}
\begin{aligned}
\left|\int_\Sigma (\mathcal H-\hbar)\Psi \ d\mu_t\right|&\leq C\left(\sigma^{-2-\frac12q-\delta q}\int_\Sigma (\mathcal H-\hbar)^2 \ d\mu+\sigma^{-1-\frac12q-\delta q}\int_\Sigma |\mathcal H-\hbar||\nabla\mathcal H| \ d\mu \right)
\\
&\leq C\sigma^{-2}\int_\Sigma (\mathcal H-\hbar)^2 \ d\mu+\varepsilon\int_\Sigma |\nabla\mathcal H|^2 \ d\mu.
\end{aligned}
\end{equation}
We conclude choosing $\varepsilon$ suitably small.\\
\indent We finally compute, using Lemma \ref{evolution41}, the evolution  
\begin{equation*}
\begin{aligned}
\frac{\textup{d}}{\textup{d}t}\int_\Sigma (\mathcal H-\hbar)^4 \ d\mu_t = & \ 4\int_\Sigma (\mathcal H-\hbar)^3\left(\Delta \mathcal H+(\mathcal H-\hbar)(|A|^2+\overline{\textup{Ric}}(\nu,\nu))-\dot\hbar\right) \ d\mu_t
\\
&+4\int_\Sigma (\mathcal H-\hbar)^3(\Upsilon-1)(\Delta\mathcal H+(\mathcal H-\hbar)(|A|^2+\overline{\textup{Ric}}(\nu,\nu)))\ d\mu_t
\\
&-\int_\Sigma H(\mathcal H-\hbar)^5 \ d\mu+\int_\Sigma \Psi(\mathcal H-\hbar)^3 \ d\mu_t.
\end{aligned}
\end{equation*}
We obtain the desired inequality \eqref{eq41645} using integration  
\begin{equation}
\int_\Sigma (\mathcal H-\hbar)^3\Delta \mathcal H \ d\mu_t=-3\int_\Sigma (\mathcal H-\hbar)^2|\nabla\mathcal H|^2 \ d\mu_t,
\end{equation}
together with the estimate
\begin{equation*}
\begin{aligned}
|\dot\hbar|\int_\Sigma |\mathcal H-\hbar|^3&\leq c\sigma^{-\frac72-\frac12q-\delta-\delta q}\left(B_\infty\sigma^{-\frac32-\delta}\right)\|\mathcal H-\hbar\|_2^2\leq cB_\infty \sigma^{-5-\frac32 q-2\delta-3\delta q},
\end{aligned}
\end{equation*}
and
\begin{align*}
\left|\int_\Sigma (\mathcal H-\hbar)^3(\Upsilon-1)\Delta\mathcal H \ d\mu\right|&\leq c\sigma^{-\frac12q-q\delta}\int_\Sigma (\mathcal H-\hbar)^2|\nabla \mathcal H|^2 \ d\mu\\
&\quad +\left|-\int_\Sigma (\mathcal H-\hbar)^3\nabla\Upsilon\nabla\mathcal H \ d\mu\right|\\
&\leq c\sigma^{-\frac12q-q\delta}\int_\Sigma (\mathcal H-\hbar)^2|\nabla \mathcal H|^2 \ d\mu\\
&\quad +c\int_\Sigma \sigma^{-\frac{q}2-1-q\delta}|\nabla\mathcal H||\mathcal H-\hbar|^3 \ d \mu.
\end{align*}
Finally, we conclude combining \eqref{estimatePsipre} and \eqref{mainpineq59}, and thus estimating
\begin{equation}
\begin{aligned}
\left|\int_\Sigma \Psi(\mathcal H-\hbar)^3 \ d\mu_t\right|&\leq C\sigma^{\frac12-\frac12q-\delta q+\delta}\int_{\Sigma}\sigma^{-\frac52-\delta}|\mathcal H-\hbar|^4+\sigma^{-\frac32-\delta}|\mathcal H-\hbar|^3|\nabla\mathcal H| \ d\mu_t\\
&\leq C\sigma^{-2}\int_\Sigma (\mathcal H-\hbar)^4 \ d\mu_t+C\sigma^{-1-\frac12q-\delta q}\int_\Sigma (\mathcal H-\hbar)^2|\mathcal H-\hbar||\nabla\mathcal H| \ d\mu_t\\
&\leq C\sigma^{-2}\int_\Sigma (\mathcal H-\hbar)^4 \ d\mu_t+C\sigma^{-q-2\delta q}\int_\Sigma (\mathcal H-\hbar)^2 |\nabla\mathcal H|^2 \ d\mu_t,
\end{aligned}
\end{equation}
where we used Young's inequality in the latter estimate. The conclusion holds for $\sigma$ suitably large.
\end{proof}
A similar estimate, but independent of the evolution of $\hbar$, can be also given for $\nabla\mathcal H$. The hypothesis on $\overset{\circ}{A}$ and the $H^1$-norm of $\mathcal H-\hbar$ are not needed in order to prove the following Lemma.
\begin{lem}\label{ineqnablaHcal} Let $(\Sigma,F_t)$, $t\in [0,T]$, such that \eqref{radius1} and \eqref{radius2} hold for every $t\in [0,T]$. Then there exists a constant $C=C(\overline c)>0$ such that 
\begin{equation}\label{eq131}
\frac{\textup d}{\textup dt}\int_\Sigma |\nabla\mathcal H|^2 \ d\mu_t\leq -\frac12\int_\Sigma |\nabla^2\mathcal H|^2 \ d\mu_t+C\sigma^{-2}\int_\Sigma |\nabla\mathcal H|^2 \ d\mu_t+C\sigma^{-4}\int_\Sigma (\mathcal H-\hbar)^2 \ d\mu_t,
\end{equation}
and
\begin{equation}\label{eq132v431}
\frac{\textup d}{\textup dt}\int_\Sigma |\nabla \mathcal H|^4 \ d\mu_t\leq -\int_\Sigma |\nabla^2\mathcal H|^2|\nabla\mathcal H|^2 \ d\mu_t+C\sigma^{-2}\int_\Sigma |\nabla \mathcal H|^4 \ d\mu_t+C\sigma^{-6}\int_\Sigma (\mathcal H-\hbar)^4 \ d\mu_t.
\end{equation}
\end{lem}
\begin{proof} We start proving inequality \eqref{eq131}. 
From Lemma \ref{stevolequations52}, after integration we get 
\begin{equation*}
\begin{aligned}
\frac{\textup d}{\textup dt}\int_\Sigma |\nabla\mathcal H|^2\ d\mu_t= & \ 2\int_\Sigma \langle\nabla\left(\Delta \mathcal H+(\mathcal H-\hbar)(|A|^2+\overline{\textup{Ric}}(\nu,\nu)\right),\nabla\mathcal H\rangle \ d\mu_t
\\
&+2\int_\Sigma \langle \nabla\left(\left(\Upsilon-1\right)\left(\Delta\mathcal H+(\mathcal H-\hbar)(|A|^2+\overline{\textup{Ric}}(\nu,\nu))\right)\right),\nabla\mathcal{H}\rangle \ d\mu_t
\\
&+\int_\Sigma |\nabla\mathcal H|^2 H(\hbar-\mathcal H) \ d\mu_t+2\int_\Sigma (\mathcal H-\hbar)|\nabla\mathcal H|^2 h^{ij}\nabla_i\mathcal H\nabla_j\mathcal H \ d\mu_t\\
&+2\int_\Sigma \langle\nabla \Psi,\nabla\mathcal H\rangle d\mu_t.
\end{aligned}
\end{equation*}
Since $H$, $|\mathcal H-\hbar|$ are bounded by $C\sigma^{-1}$ and $\left||A|^2+\overline{\textup{Ric}}(\nu,\nu)\right|\leq C\sigma^{-2}$ and $|\textup{Ric}^\Sigma|\leq C\sigma^{-2}$, using Bochner's identity and integration by part we get 
\begin{equation}
\begin{aligned}
\frac{\textup d}{\textup dt}\int_\Sigma |\nabla\mathcal H|^2\ d\mu_t\leq & \ C \sigma^{-2}\int_\Sigma |\nabla\mathcal H|^2 \ d\mu_t+C\sigma^{-4}\int_\Sigma (\mathcal H-\hbar)^2 \ d\mu_t-\int_\Sigma |\nabla^2\mathcal H|^2 \ d\mu_t
\\
&+2\left|\int_\Sigma \left(\Upsilon-1\right)\left(\Delta\mathcal H+(\mathcal H-\hbar)(|A|^2+\overline{\textup{Ric}}(\nu,\nu))\right)\Delta\mathcal H \ d\mu_t\right|
\\
&+2\left|\int_\Sigma \Psi\Delta \mathcal H \ d\mu_t\right|.
\end{aligned}
\end{equation}
Note that, combining \eqref{estimatePsipre} and \eqref{mainpineq59}, we get
\begin{equation}
\begin{aligned}
\left|\int_\Sigma \Psi\Delta\mathcal H \ d\mu_t\right|\leq & \ C\sigma^{\frac12-\frac12q-\delta q+\delta}\int_\Sigma\left(\sigma^{-\frac52-\delta}|\mathcal H-\hbar|+\sigma^{-\frac32-\delta}|\nabla\mathcal H|\right)|\nabla^2\mathcal H| \ d\mu_t
\\
\leq & \ C\sigma^{-4}\int_\Sigma (\mathcal H-\hbar)^2 \ d\mu_t+C\sigma^{-2}\int_\Sigma|\nabla\mathcal H|^2 \ d\mu_t+C\sigma^{-q-2\delta q}\int_\Sigma |\nabla^2\mathcal H|^2 \ d\mu_t.
\end{aligned}
\end{equation}
Since $|\Upsilon-1|=\textup O(\sigma^{-\frac12q-q\delta})$  and using Young's inequality, we conclude, for $\sigma$ large,
\begin{equation}
\begin{aligned}
\frac{\textup d}{\textup dt}\int_\Sigma |\nabla\mathcal H|^2\ d\mu_t\leq & \ C \sigma^{-2}\int_\Sigma |\nabla\mathcal H|^2 \ d\mu_t+C\sigma^{-4}\int_\Sigma (\mathcal H-\hbar)^2 \ d\mu_t-\int_\Sigma |\nabla^2\mathcal H|^2 \ d\mu_t
\\
&+C\sigma^{-\frac12q-q\delta}\int_\Sigma |\nabla^2\mathcal H|^2 \ d\mu_t+C\sigma^{-2-\frac12q-q\delta}\int_\Sigma |\mathcal H-\hbar||\nabla^2\mathcal H| \ d\mu_t
\\
&+C\sigma^{-q-2q\delta}\int_\Sigma |\nabla^2\mathcal H|^2 \ d\mu_t\\
\leq & \ C \sigma^{-2}\int_\Sigma |\nabla\mathcal H|^2 \ d\mu_t+C\sigma^{-4}\int_\Sigma (\mathcal H-\hbar)^2 \ d\mu_t-\frac12\int_\Sigma |\nabla^2\mathcal H|^2 \ d\mu_t.
\end{aligned}
\end{equation}
\indent We now prove \eqref{eq132v431}. From \eqref{eq101nbis} we get, after integrating by parts, 
\begin{equation}\label{ev72}
\begin{aligned}
\frac{\textup d}{\textup d t}\int_\Sigma |\nabla\mathcal H|^4 \ d\mu_t &= 4\int_\Sigma (\mathcal H-\hbar)|\nabla\mathcal H|^2h^{ij}\nabla_i\mathcal H\nabla_j\mathcal H\\
&\quad +4\int_\Sigma \langle\nabla\left(\Delta \mathcal H+(\mathcal H-\hbar)(|A|^2+\overline{\textup{Ric}}(\nu,\nu)\right),\nabla\mathcal H\rangle \ d\mu_t\\
&\quad +\int_\Sigma |\nabla \mathcal H|^4 H(\hbar-\mathcal H) \ d\mu_t
\\
&\quad -4\int_\Sigma \left(\Upsilon-1\right)\left(\Delta \mathcal H+(\mathcal H-\hbar)(|A|^2+\overline{\textup{Ric}}(\nu,\nu)\right)\Delta\mathcal H|\nabla\mathcal H|^2 \ d\mu_t\\
&\quad +4\int_\Sigma \langle \nabla\Psi,\nabla\mathcal H\rangle \ d\mu_t.
\end{aligned}
\end{equation}
To estimates the terms above, note that, if $\sigma_0$ is so large that $|\Upsilon-1|\leq \varepsilon$ (see \eqref{Upsilonineq4545}), then 
\begin{align*}
&\left|4\int_\Sigma \left(\Upsilon-1\right)\left(\Delta \mathcal H+(\mathcal H-\hbar)(|A|^2+\overline{\textup{Ric}}(\nu,\nu)\right)\Delta\mathcal H|\nabla\mathcal H|^2 \ d\mu_t\right|\\
&\qquad\qquad\qquad\leq \varepsilon\int_\Sigma |\nabla^2\mathcal H|^2|\nabla \mathcal H|^2 \ d\mu_t +C\sigma^{-2}\int_\Sigma |\mathcal H-\hbar||\nabla^2 \mathcal H||\nabla\mathcal H|^2 \ d\mu_t.
\end{align*}
Moreover, using again integration by parts on $\langle\nabla\Psi,\nabla\mathcal H\rangle=\nabla\cdot\left(\Psi\nabla\mathcal H\right)-\Psi\Delta\mathcal H$, we estimate 
\begin{equation*}
\begin{aligned}
\left|\int_\Sigma \langle \nabla\Psi,\nabla\mathcal H\rangle|\nabla\mathcal H|^2 \ d\mu_t\right|&\leq C\sigma^{-1-\frac12q-\delta q}\left(\sigma^{-1}\int_\Sigma |\mathcal H-\hbar||\nabla^2\mathcal H||\nabla \mathcal H|^2 \ d\mu_t\right)\\
&\quad +C\sigma^{-1-\frac12q-\delta q}\int_\Sigma |\nabla^2\mathcal H||\nabla\mathcal H|^3 \ d\mu_t,
\end{aligned}
\end{equation*}
where we also used \eqref{estimatePsipre} combined with \eqref{mainpineq59}. We conclude using Bochner's formula, the inequality $H+|H-h|+|A|\leq C\sigma^{-1}$ and that $|\textup{Ric}^\Sigma|\leq C\sigma^{-2}$, obtaining 
\begin{equation}
\begin{aligned}
\frac{\textup d}{\textup d t}\int_\Sigma |\nabla\mathcal H|^4 \ d\mu_t&\leq-4\int_\Sigma |\nabla^2\mathcal H|^2|\nabla \mathcal H|^2+C\sigma^{-2}\int_\Sigma |\nabla \mathcal H|^4 \ d\mu
\\
&\quad +C\sigma^{-2}\int_\Sigma |\mathcal H-\hbar||\nabla^2 \mathcal H||\nabla\mathcal H|^2 \ d\mu_t+\varepsilon\int_\Sigma |\nabla^2\mathcal H|^2|\nabla \mathcal H|^2 \ d\mu_t\\
&\quad +C\sigma^{-1}\int_\Sigma |\nabla^2\mathcal H||\nabla\mathcal H|^3 \ d\mu_t.
\end{aligned}
\end{equation}
The desired inequality appears when using Young's inequality
\begin{flushleft}
$\displaystyle C\sigma^{-2}\int_\Sigma |\mathcal H-\hbar||\nabla^2 \mathcal H||\nabla\mathcal H|^2 \ d\mu_t+C\sigma^{-1}\int_\Sigma |\nabla^2\mathcal H||\nabla\mathcal H|^3 \ d\mu_t$
\end{flushleft}
\begin{equation*}
\qquad \qquad \quad \leq \varepsilon\int_\Sigma |\nabla^2 \mathcal H|^2|\nabla\mathcal H|^2 \ d\mu_t+C\sigma^{-6}\int_\Sigma (\mathcal H-\hbar)^4 \ d\mu_t+C\sigma^{-2}\int_\Sigma |\nabla\mathcal H|^4 \ d\mu_t
\end{equation*}
and choosing $\varepsilon$ suitably small.
\end{proof}
\begin{lem}\label{useful2} Let $\Sigma\hookrightarrow \M$ be a \textup{surface}. Then we have, for every $\varepsilon>0$ and $\sigma>1$,  
\begin{equation*}
-\sigma^{-4}\int_\Sigma (\mathcal H-\hbar)^2|\nabla \mathcal H|^2 \ d\mu\leq -\frac{\varepsilon}{2\sigma^2}\int_\Sigma |\nabla \mathcal H|^4 \ d\mu+\varepsilon^2\int_\Sigma |\nabla^2\mathcal H|^2|\nabla \mathcal H|^2 \ d\mu.
\end{equation*}
\end{lem} 
The proof of Lemma \ref{useful2} is an easy consequence of Young's inequality and integration by parts. See \cite{vpmcf} for details. This leads to the following Lemma.
\begin{lem}\label{lemtech9}
Let $(\Sigma,F_t)$, $t\in[0,T]$, such that \eqref{radius1}, \eqref{radius2}, \eqref{Binftycinftyst}, \eqref{H^1est411} and $\|\overset{\circ}{A}\|_{L^4(\Sigma_t)}\leq B_1 \sigma^{-1-\delta}$ for every $t\in [0,T]$. For $\eta>0$, let us set 
\begin{equation}
a_\eta(t):=\eta\sigma^{-4}\|\mathcal H-\hbar\|_{L^4(\Sigma_t)}^4+\|\nabla\mathcal H\|_{L^4(\Sigma_t)}^4.
\end{equation}
Then there exist a universal constant $\eta_w>0$ and a radius $\sigma_0=\sigma_0(B_1,\cin,\overline c,\delta)>1$ such that if $\eta=\eta_w$ and $\sigma>\sigma_0$ the following statements hold.
\begin{enumerate}[label=\textup{(\roman*)}]
\item  There exists a constant $c=c(B_1,\eta_w,\overline c)$ such that if $B_2>c$ we have the implication 
\begin{equation*}
a_{\eta_w}(0)<B_2\sigma^{-8-4\delta} \ \Longrightarrow \ a_{\eta_w}(t)<B_2\sigma^{-8-4\delta} \ \textup{for every $t\in [0,T]$}.
\end{equation*}
\item  If in addition we suppose \eqref{Binftycinftyst}, there exists a constant $c=c(\cin,B_\infty)$ such that if we choose $B_{\textup{in}}>c(\cin,B_\infty)$ we have the implication
\begin{equation*}
\textup a_\eta(0)<B_{\textup{in}} \sigma^{-7-\frac32q-2\delta-3\delta q} \ \Longrightarrow \  \textup a_\eta(t)<B_{\textup{in}} \sigma^{-7-\frac32q-2\delta-3\delta q} \ \textup{for every $t\in[0,T]$}.
\end{equation*}
\end{enumerate}
\end{lem}
\begin{proof} Combining Lemma \ref{rateh(t)st}, Lemma \ref{ineqnablaHcal} and Lemma \ref{useful2}, we get 
\begin{equation}\label{eq182}
\dot a_\eta(t)\leq -C\sigma^{-2}a_\eta(t)+\tilde C\sigma^{-6}\int_\Sigma (\mathcal H-\hbar)^4 \ d\mu_t+\tilde cB_\infty \sigma^{-9-\frac32 q-2\delta-3\delta q},
\end{equation}
for some $\tilde C$ universal constant and $\tilde c=\tilde c(\cin,\overline c)$. We will use inequality \eqref{eq182} in order to prove two different conclusions.
\begin{enumerate}[label=\textup{(\roman*)}]
\item Since $q\geq 2$, choosing $\sigma$ suitably large depending on $B_\infty$, such that it holds $\tilde cB_\infty \sigma^{-9-\frac32 q-2\delta-3\delta q}\leq \sigma^{-10-4\delta}$, we have 
\begin{equation}\label{eqa2.0}
\dot a_\eta(t)\leq -C\sigma^{-2}a_\eta(t)+\tilde C\sigma^{-6}\int_\Sigma (\mathcal H-\hbar)^4 \ d\mu_t+\sigma^{-10-4\delta}.
\end{equation}
Moreover, Lemma \ref{cor1} implies that
\begin{equation*}
\int_\Sigma (H-h)^4 \ d\mu_t\leq \cper^4\left(\|\overset{\circ}{A}\|_{L^4(\Sigma,\mu_t)}^4+\sigma^{-4-4\delta}\right)\leq \cper^4(B_1^4+1)\sigma^{-4-4\delta}.
\end{equation*} 
and thus \eqref{eqa2.0} becomes
\begin{equation}\label{eq154bis}
\dot a_\eta(t)\leq -C\sigma^{-2}a_\eta(t)+c\sigma^{-10-4\delta}
\end{equation}
with $c=c(B_1,\overline c,\cper)$. Thus, if $B_2>c/C$, we get the thesis. 
\item Since we are assuming \eqref{H^1est411} for every $t\in[0,T]$, the Sobolev's immersion (see \cite[Lemma 12]{cederbaum}) implies that 
\begin{equation*}
\|\mathcal H-\hbar\|_{L^4(\Sigma_t)}\leq \tilde c\sigma^{-\frac12-\frac{q}{2}-q\delta} \quad \textup{for every $t\in [0,T]$},
\end{equation*}
where $\tilde c=\tilde c(\cin,\overline c)$. Thus \eqref{eq182} becomes 
\begin{equation}
\begin{aligned}
\dot a_\eta(t)\leq& -C\sigma^{-2}a_\eta(t)+\tilde c\sigma^{-8-2q-4q\delta}+\tilde cB_\infty \sigma^{-9-\frac32 q-2\delta-3\delta q}\\
\leq& -C\sigma^{-2}a_\eta(t)+2\tilde cB_\infty \sigma^{-9-\frac32q-2\delta-3\delta q}
\end{aligned}
\end{equation}
for $\sigma$ large, since $q\geq 2$ and $\delta\in (0,\frac12]$. Choosing $B_\textup{in}>2\tilde c B_\infty/C$ we have the thesis.
\end{enumerate}
\end{proof}
From now on, when considering the roundness class $\mathcal W_\sigma^\eta(B_1,B_2)$, we fix the parameter $\eta$ equal
to the value $\eta_w$ given by the previous Lemma, and we will no longer need to specify the dependence
on $\eta$ of the constants in the estimates.
\subsection{Evolution of $\|\mathcal H-\hbar\|_{L^2(\Sigma_t)}$ and convergence} An important assumption in the previous results was the comparability between $r_\Sigma$ and $\sigma$ in \eqref{radius2} which assures that on $\Sigma_t$ the ambient curvature decays with the right order, as highlighted in Remark \ref{preliminaryremark27}. To justify this assumption, we study now the evolution of $L^2$-norm of $\mathcal H-\hbar$, which relies on the spectral analysis of Section 3. 
The following Lemma is an improvement of inequality \eqref{roughestL2t}. Under an additional hypothesis, this inequality shows that the negative term in the evolution of $\|\mathcal H-\hbar\|_{L^2(\Sigma_t)}^2$ is dominant.
\begin{lem}\label{gronwallst} Let $(\Sigma,F_t)$, $t\in [0,T]$, be such that \eqref{radius1}, \eqref{radius2}, \eqref{Binftycinftyst}, \eqref{H^1est411} hold for every $t\in [0,T]$. For every $\Omega>0$ there exists $\sigma_0(\overline c,\Omega)>1$ such that if 
\begin{equation}\label{hypomega558}
\|\mathcal H-\hbar\|_{L^\infty(\Sigma_t)}\leq \Omega\sigma^{-\frac54-\frac38q-\frac{\delta}2-\frac{3\delta q}{4}}, \quad \forall t\in [0,T]
\end{equation}
and $\sigma>\sigma_0$, then 
\begin{equation*}
\frac{\textup d}{\textup dt}\int_\Sigma (\mathcal H-\hbar)^2 \ d\mu_t\leq -\frac{E_\textup{ADM}}{\sigma^3}\int_\Sigma (\mathcal H-\hbar)^2 d \mu_t,
\end{equation*}
for every $t\in [0,T]$.
\end{lem}
\begin{proof} We easily compute the evolution
\begin{equation}\label{maniL2559}
\begin{aligned}
\frac{\textup d}{\textup dt}\int_\Sigma (\mathcal H-\hbar)^2 \ d\mu_t =&  -2\int_\Sigma (\mathcal H-\hbar)L\left(\mathcal H-\hbar\right) \ d\mu_t
\\
&+2\int_\Sigma (\mathcal H-\hbar)(\Upsilon-1)\left(\Delta\mathcal H+(\mathcal H-\hbar)(|A|^2+\overline{\textup{Ric}}(\nu,\nu))\right) \ d\mu_t
\\
&+\int_\Sigma \Psi(\mathcal H-\hbar) \ d\mu_t-\int_\Sigma H(\mathcal H-\hbar)^3 \ d\mu_t.
\end{aligned}
\end{equation}
Combining \eqref{estimatePsipre} and \eqref{mainpineq59} we get 
\begin{equation}
\left|\Psi\right|\leq C(\sigma|P|)^{q-1}\left(\sigma^{-\frac52-\delta}|\mathcal H-\hbar|+\sigma^{-\frac32-\delta}|\nabla\mathcal H|\right),
\end{equation}
which implies, using that $\sigma|P|\leq C\sigma^{-\frac12-\delta}$,
\begin{equation*}
\begin{aligned}
\left|\int_\Sigma (\mathcal H-\hbar)\Psi \ d\mu_t\right|&\leq C\sigma^{-2-\frac12q-\delta q}\int_\Sigma (\mathcal H-\hbar)^2 \ d\mu_t+C\sigma^{-\frac32-\delta}\int_\Sigma (\sigma|P|)^{q-1}|\mathcal H-\hbar||\nabla\mathcal H| \ d\mu_t \\
&\leq \frac{\varepsilon E_\textup{ADM}}{\sigma^3}\int_\Sigma (\mathcal H-\hbar)^2 \ d\mu_t+\varepsilon\int_\Sigma (\sigma|P|)^{2q-2}|\nabla\mathcal H|^2 \ d\mu_t.
\end{aligned}
\end{equation*}
where in the latter inequality we used parametric Young's inequality and we have chosen $\sigma$ suitably large.\\
\indent We now estimate, using integration by parts and formula \eqref{eq561p} for $\nabla\Upsilon$, together with the fact that $\left|\frac{H}{\mathcal H}\right|\leq C$ and $H\sim \mathcal H\sim \frac2\sigma$,
\begin{equation}\label{eq5611}
\begin{aligned}
 \int_\Sigma (\mathcal H-\hbar)(\Upsilon-1)\Delta\mathcal H \ d\mu_t &\leq -\int_\Sigma (\Upsilon-1)|\nabla\mathcal H|^2 \ d\mu_t\\
&\quad +\int_\Sigma \left(\sigma^{q+1}|P|^q|\nabla\mathcal H|+\sigma^q|P|^{q-1}|\nabla P|\right)|\mathcal H-\hbar||\nabla\mathcal H| \ d\mu_t.
\end{aligned}
\end{equation}
Since $|\mathcal H-\hbar|\leq \epsilon\sigma^{-1}$ for $\sigma$ suitably large, and $|\nabla P|\leq \sigma^{-\frac52-\delta}$ because of Lemma \ref{propertiesinitialdata},
\begin{equation}\label{eq5612}
\begin{aligned}
&\left(\sigma^{q+1}|P|^q|\nabla\mathcal H|+\sigma^q|P|^{q-1}|\nabla P|\right)|\mathcal H-\hbar||\nabla\mathcal H|
\\
&\qquad \qquad\leq \epsilon \sigma^q|P|^q|\nabla \mathcal H|^2+\sigma^{q-\frac52-\delta}|P|^{q-1}|\mathcal H-\hbar||\nabla\mathcal H|\\
&\qquad \qquad\leq \epsilon(\sigma|P|)^q|\nabla\mathcal H|^2+\sigma^{-\frac32-\delta}(\sigma|P|)^{q-1}|\mathcal H-\hbar||\nabla\mathcal H|\\
&\qquad \qquad\leq \epsilon(\sigma|P|)^q|\nabla\mathcal H|^2+C\sigma^{-3-2\delta}|\mathcal H-\hbar|^2+\epsilon (\sigma|P|)^{2q-2}|\nabla\mathcal H|^2
\end{aligned}
\end{equation}
where in the latter inequality we used parametric Young's inequality. Since it holds $(\sigma|P|)^{2q-2}=(\sigma|P|)^{q-2}\left(\sigma|P|\right)^q\leq C(\sigma|P|)^q$, combining \eqref{eq5611} and \eqref{eq5612} we get
\begin{equation}
\begin{aligned}
\int_\Sigma (\mathcal H-\hbar)(\Upsilon-1)\Delta\mathcal H \ d\mu_t &\leq -\int_\Sigma (\Upsilon-1)||\nabla\mathcal H|^2 \ d\mu_t\\
&\quad +C\sigma^{-3-2\delta}\int_\Sigma (\mathcal H-\hbar)^2 \ d\mu_t+\epsilon C\int_\Sigma (\sigma|P|)^q|\nabla\mathcal H|^2 \ d\mu_t,
\end{aligned}
\end{equation}
Note furthermore that \eqref{hypomega558} implies
\begin{equation*}
\left|\int_\Sigma H(\mathcal H-\hbar)^3\ d\mu_t\right|\leq
 \Omega\sigma^{-\frac94-\frac38q-\frac{3q\delta}4-\frac\delta2}\int_\Sigma (\mathcal H-\hbar)^2 \ d\mu_t\leq \frac{\varepsilon E}{\sigma^3}\int_\Sigma (\mathcal H-\hbar)^2 \ d\mu_t,
\end{equation*}
if $\sigma>\sigma_0$, for some $\sigma_0=\sigma_0(\Omega)$.\\
\indent We conclude from \eqref{maniL2559}, using Proposition \ref{spacetimespectrath} and $|\Upsilon-1||A|^2\leq C\sigma^{-2-\frac{q}2-q\delta}\leq \varepsilon E_\textup{ADM}\sigma^{-3}$ and $\Upsilon-1\geq \underline c\sigma^q|P|^q$ because of \eqref{Upsilonineq4545}, obtaining
\begin{equation}
\begin{aligned}
\frac{\textup d}{\textup dt}\int_\Sigma (\mathcal H-\hbar)^2 \ d\mu_t &\leq -\frac{2E_\textup{ADM}}{\sigma^3}\int_\Sigma (\mathcal H-\hbar)^2 \ d\mu_t+\frac{\varepsilon E_\textup{ADM}}{\sigma^3}\int_\Sigma (\mathcal H-\hbar)^2 \ d\mu_t\\
& \quad \ +(\epsilon C-\underline c)\int_\Sigma \left(\sigma|P|\right)^q|\nabla\mathcal H|^2 \ d\mu_t+C\sigma^{-3-2\delta}\int_\Sigma (\mathcal H-\hbar)^2 \ d\mu_t\\
&\leq -\frac{E_\textup{ADM}}{\sigma^3}\int_\Sigma (\mathcal H-\hbar)^2 \ d\mu_t
\end{aligned}
\end{equation}
for $\epsilon$ small with respect to $\underline c$, $\varepsilon<1$ and $\sigma$ suitably large.
\end{proof}
The next result, which is similar to Proposition 3.4 in \cite{huiskenyau}, gives a bound on the possible change
of area of the surface along the flow as long as it remains round. The proof is analogous to that of \cite[Lemma 4.12]{vpmcf}.
\begin{lem}\label{errorradius44}
Given $B_1,B_2$, there exist constants $c>0$ and $\sigma_0>1$ such that, if $\sigma>\sigma_0$ and $\Sigma_t$ is a solution of the flow \eqref{generalflow513} for $t \in [0,T]$ with $\Sigma_t \in  \mathcal{W}^\eta_\sigma(B_1,B_2)$ for all $t$ then
\begin{equation*}
|\sigma_{\Sigma_0}-\sigma_{\Sigma_t}|\leq c\sigma^{\frac12-\delta}
\end{equation*}
for every $t\in[0,T]$.
\end{lem}
We are now ready to prove that, by an appropriate choice of the parameters of roundness class, a well-centered round surface remains inside the class for arbitrary times. Remember that hypotheses \eqref{hypstmainth2in} are in particular satisfied by Nerz's foliation.\\
\indent The following Theorem is the key step in the proof of Theorem \ref{mainthm}.
\begin{thm}\label{proofmainth411} Let $(\M,\overline\g,\overline\K)$ be a $C_{\frac12+\delta}^2$-asymptotically flat initial data set, with $E_\textup{ADM}>0$. Choose $B_1$ as in Lemma \ref{eqevApallST} and $B_2$ and $\eta$ as in Lemma \ref{lemtech9}. For every $C_0>0$ there exist $\overline B=\overline B(C_0)$ and $\sigma_0=\sigma_0(\overline c,\delta,E_\textup{ADM},B_1,B_2,C_0)$ such that the following holds. Let $(\Sigma,F_t)$ be a solution to the volume preserving spacetime mean curvature flow for $t\in [0,T]$ such that $\Sigma_0$ \textup{(i)} belongs to $\mathcal B_\sigma(B_1,B_2,B_\textup{cen})$ with $\sigma=\sigma_{\Sigma_0}$, \textup{(ii)} is a CMC-surface and \textup{(iii)} $|\vec z_{\Sigma_0}|\leq C_0\sigma^{1-\delta}$. Then, if $B_\textup{cen}\geq \overline B$ and $\sigma\geq\sigma_0$, $\Sigma_t$ belongs to $\mathcal B_\sigma(B_1,B_2,B_\textup{cen})$ for every $t\in [0,T]$.
\end{thm}
\begin{rem}
Note that the following proof also works when $\Sigma$ is almost CMC and not exactly CMC.
\end{rem}
\begin{proof} Note that, since $\Sigma_0=\Sigma$ belongs to $\mathcal B_\sigma(B_1,B_2,B_\textup{cen})$, then it satisfies 
\begin{equation*}
\frac12\leq \frac{r_\Sigma}{\sigma}\leq 2, \qquad 1\leq \sigma H\leq \sqrt 5.
\end{equation*}
Thus Lemma \ref{propertiesinitialdata} implies that the initial (CMC) surface satisfies 
\begin{equation}\label{H1hypC1}
\|\mathcal H-\hbar\|_{H^1(\Sigma_0)}\leq C_1 \sigma^{-\frac12q-q\delta}
\end{equation}
for some $C_1=C_1(\overline c)>0$. Thus, we define $T_\textup{max}$ as the supremum of the times $\overline T\leq T$ such that the following conditions hold:
\begin{enumerate}[label=(\roman*)]
\item $F_t$ exists in $[0,\overline T]$;
\item $\|\mathcal H-\hbar\|_{L^2(\Sigma_t)}<(C_1+1)\sigma^{-\frac12q-q\delta}$ for every $t\in [0,\overline T]$;
\item $\Sigma_t \in \mathcal B_{\sigma}(B_1,B_2,B_\textup{cen})$ for every $t\in [0,\overline T]$.
\end{enumerate}
We remark that $T_\textup{max}>0$ and $\Sigma_{T_\textup{max}}$ belongs to $\overline{\mathcal B}_\sigma(B_1,B_2,B_\textup{cen})\subset \overline{\mathcal W}_\sigma(B_1,B_2)$. By Lemma \ref{cor1}, Lemma \ref{errorradius44} and Definition \ref{roundsurface}, the conditions \eqref{radius1} and \eqref{radius2} hold for every $t\in[0,T_\textup{max}]$. See again Remark \eqref{remark26} for a direct estimate of the Euclidean radius.\\
\\
\textbf{Claim:} There exists $\cin=\cin(\overline c)>0$ such that \eqref{H^1est411} holds for every $t\in [0,T_\textup{max}]$.\\
\noindent \textit{Proof of the Claim.} Combining together \eqref{roughestL2t} and \eqref{eq131}, if $C$ is the maximum between the two constants involved, we find that 
\begin{equation}
\begin{aligned}
\frac{\textup d}{\textup dt}\left(\int_\Sigma |\nabla\mathcal H|^2 \ d\mu_t+4C\sigma^{-2}\int_\Sigma (\mathcal H-\hbar)^2 \ d\mu_t\right)\leq & -\frac12\int_\Sigma |\nabla^2\mathcal H|^2 \ d\mu_t-C\sigma^{-2}\int_\Sigma |\nabla\mathcal H|^2 \ d\mu_t
\\
&+(4C^2+C)\sigma^{-4}\int_\Sigma (\mathcal H-\hbar)^2 \ d\mu_t.
\end{aligned}
\end{equation}
Setting $a(t):=\|\nabla \mathcal H\|_{L^2(\Sigma_t)}^2+4C\sigma^{-2}\|\mathcal H-\hbar\|_{L^2(\Sigma_t)}^2$, since by definition of $T_\textup{max}$ it holds $\|\mathcal H-\hbar\|_{L^2(\Sigma_t)}<(C_1+1)\sigma^{-\frac12q-q\delta}$ for every $t\in [0,T_\textup{max}]$,
\begin{equation}\label{diffeqabis}
\begin{aligned}
\dot a(t)&\leq-C\sigma^{-2} a(t)+(8C^2+C)\sigma^{-4}\int_\Sigma (\mathcal H-\hbar)^2 \ d\mu_t
\\
&\leq-C\sigma^{-2}a(t)+2(8C^2+C)(C_1+1)^2\sigma^{-4-q-2q\delta}
\end{aligned}
\end{equation}
Since, by \eqref{H1hypC1}, $a(0)\leq (1+4C)C_1\sigma^{-2-q-2\delta q}$, \eqref{diffeqabis} implies that $a(t)\leq C(\overline c,C_1)\sigma^{-2-q-2\delta q}$ for every $t\in[0,T_\textup{max}]$. Since also $C_1=C_1(\overline c)$, this proves that there exists $\cin=\cin(\overline c)$ such that the claim holds.\\
\\
\noindent Now, \eqref{radius1}, \eqref{radius2}, \eqref{Binftycinftyst} and \eqref{H^1est411} imply that we are in the hypotheses of Proposition \ref{eqevApallST} and of point (i) of Lemma \ref{lemtech9}. Thus, the choices of $B_1$ and $B_2$ imply that $\Sigma_{T_\textup{max}}$ belongs to $\mathcal W_\sigma(B_1,B_2)$ for $\sigma$ large. Moreover, \eqref{Binftycinftyst} holds for some $B_\infty=B_\infty(B_1,B_2)$, thanks again to Lemma \ref{cor1}.
\\
\indent We conclude showing that, if $B_\textup{cen}$ is chosen suitably large, then $\Sigma_t\in\mathcal B_\sigma(B_1,B_2,B_\textup{cen})$ for every $t\in [0,T_\textup{max}]$. Since $\Sigma_0$ is a CMC-surface, it is easy to verify that $a_{\eta_w}(0)<B_\textup{in}\sigma^{-7-\frac32q-2\delta-3\delta q}$ for a constant $B_\textup{in}$ suitably large. We remember that the function $a_\eta$ has been defined in Lemma \ref{lemtech9}. Moreover, Lemma \ref{lemtech9}, point (ii), implies that if $B_\textup{in}$ is chosen suitably large, depending on $\cin$ and $B_\infty$, then $\|\mathcal H-\hbar\|_{L^\infty(\Sigma_t)}\leq B_\textup{in}\sigma^{-\frac54-\frac38q-\frac\delta2-\frac{3\delta q}{4}}$ holds for every $t\in [0,T_\textup{max}]$. Thus Lemma \ref{gronwallst}, with $\Omega:=B_\textup{in}$, combined with Gronwall's Lemma, implies that 
\begin{equation}\label{expdecay457}
\|\mathcal H-\hbar\|_{L^2(\Sigma_t)}\leq \|\mathcal{H}-\hbar\|_{L^2(\Sigma_0)}e^{-\frac{E_\textup{ADM}t}{2\sigma^3}}<(C_1+1)\sigma^{-\frac12q-q\delta}e^{-\frac{E_\textup{ADM}t}{2\sigma^3}},
\end{equation}
for every $t\in [0,T_\textup{max}]$. Setting, $\vec z(t)=\vec z_{\Sigma_t}$, we show that the behavior of the barycenter is controlled. Analogously to \cite{corvinowu}, we have the evolution 
\begin{equation}\label{barycenterst449}
\partial_t\left(|\Sigma_t|\vec z(t)\right)=\int_\Sigma (\hbar-\mathcal H)\nu \ d\mu_t+\int_\Sigma F_t(x) H(\hbar-\mathcal H) \ d\mu_t.
\end{equation}
Combining this with the estimates $H\leq \frac{5}{\sigma}$, $|F_t(x)|\leq R_\Sigma(t)\leq 3\sigma$ and \eqref{expdecay457}, we obtain
\begin{equation}\label{ineq448}
\partial_t\left(|\Sigma_t||\vec z(t)|\right)\leq C\sigma\|\mathcal H-\hbar\|_{L^2(\Sigma,\mu_t)}<C(C_1+1) \sigma^{1-\frac{q}2-q\delta} e^{-\frac{E_\textup{ADM}t}{2\sigma^3}}.
\end{equation}
Integrating \eqref{ineq448} over $[0,T_\textup{max}]$, we get 
\begin{equation*}
|\Sigma_{T_\textup{max}}||\vec z(T_\textup{max})|-|\Sigma_0||\vec z(0)|\leq C(C_1+1)\sigma^{1-\frac{q}2-q\delta}\left(\frac{2\sigma^3}{E_\textup{ADM}}\right)\left(1-e^{-\frac{E_\textup{ADM}T_\textup{max}}{2\sigma^3}}\right).
\end{equation*}
By the hypotheses $|\vec z_{\Sigma_0}|\leq C_0\sigma^{1-\delta}$, we find 
\begin{equation}
|\vec z(T_\textup{max})|\leq \frac{2}{7\pi}\left(5\pi(C_0\sigma^{1-\delta})+\frac{2C(C_1+1)}{E_\textup{ADM}}\sigma^{2-\frac{q}2-q\delta}\right)<B_\textup{cen} \sigma^{1-\delta}
\end{equation}
if $B_\textup{cen}$ suitably large, depending on $C_0$, $C$, $C_1$ and $E_\textup{ADM}$. Thus $\Sigma_{T_\textup{max}}$ belongs to  the class $\mathcal B_\sigma(B_1,B_2,B_\textup{cen})$, and combining this with \eqref{expdecay457} we obtain that necessarily $T_\textup{max}=T$.
\end{proof}
\noindent\textbf{Local regularity of the flow.} We now review the regularity theory of the non-linear flow we are considering. Since in a local interval of existence $[0,t_0)$ the principal curvatures are uniformly bounded (by the preservation of the roundness), it follows that $\Sigma_t$ can be locally written as a graph. Suppose in particular that $\Sigma_t\cap B_\epsilon (x_0)=\{\left(x_1,x_2,u(t,x_1,x_2)\right): \ (x_1,x_2)\in\mathcal A\}$, with $\mathcal A\subset \R^2$ open. Since the metric, the unit normal vector, and the mean curvature of $\Sigma_t$ are locally given by
\begin{equation}
\begin{aligned}
g_{ij}=\delta_{ij}+D_iuD_ju, \qquad \nu=\frac{\left(-D_1u,-D_2u,1\right)}{\sqrt{1+|Du|^2}},\\
H=\frac{1}{\sqrt{1+|Du|^2}}\left(\delta^{ij}-\frac{D^iuD^ju}{1+|Du|^2}\right)D_{ij}^2 u,
\end{aligned}
\end{equation}
the equation \eqref{stflow1}, written in a tangential fashion, translates into an equation for $u$ 
\begin{equation}\label{equationforu}
\partial_t u=\sqrt{1+|Du|^2}\left(\Phi\left(\frac{1}{\sqrt{1+|Du|^2}}\left(\delta^{ij}-\frac{D^iuD^ju}{1+|Du|^2}\right)D_{ij}^2 u,P\right)-\hbar\right),
\end{equation}
where $P=g^{ij}\overline\K_{ij}$ is a smooth function and $\Phi(s,\gamma)=\sqrt[q]{s^q-|\gamma|^q}$. We rewrite equation \eqref{equationforu} as
\begin{equation}
\partial_t u=\mathcal F(D^2u,Du,x,t).
\end{equation}
Note that 
\begin{equation}
\dot{\mathcal F}^{ij}:=\frac{\partial \mathcal F}{\partial D_{ij}^2 u}=\left(\delta^{ij}-\frac{D^iuD^ju}{1+|Du|^2}\right)\partial_s \Phi
\end{equation}
and $\partial_s\Phi=q^{-1}\left(s^q-|\gamma|^q\right)^{\frac1q-1}(qs^{q-1})>0$. Thus, as a matrix, 
\begin{equation}
\begin{aligned}
|w|^2\left(\inf_{\mathcal A}\partial_s \Phi\right)\left(1-\frac{\sup|Du|^2}{1+\sup|Du|^2}\right)&\leq\dot{\mathcal F}^{ij}w_iw_j=\left(|w|^2-\frac{(Du\cdot w)^2}{1+|Du|^2}\right)\partial_s\Phi\\
&\leq |w|^2\left(\sup_{\mathcal A}\partial_s \Phi\right).
\end{aligned}
\end{equation}
Finally note that if $\dot {\mathcal F}^{ij}M_{ij}=0$ then $\left(\delta^{ij}-\frac{D^iuD^ju}{1+|Du|^2}\right)M_{ij}=0$. Thus, computing $\ddot{\mathcal F}^{ij,kl}:=\frac{\partial^2 \mathcal F}{\partial D_{ij}^2 u\partial D_{kl}^2 u}$, it follows that this implies that $\ddot{\mathcal F}^{ij,kl}M_{ij}M_{kl}=0$.\\
\indent This means that we are in the hypothesis of Theorem 6 in \cite{benandrews}, which let us obtain a $C^{2,\alpha}$ estimate on $u$, for some $\alpha\in(0,1)$. By standard arguments, this means that the coefficients of the linearization of the non-linear equation are $C^{0,\alpha}$-H\"older, and thus the standard theory (see for example \cite{ladysoloural}) implies uniform bounds on all higher derivatives of $u$. Covering $\Sigma_t$ with graphs over balls of the same radius, we obtain H\"older estimates on the curvature and its derivatives.\\
\textbf{Proof of Theorem \ref{mainthm}.} Consider a CMC surface $\Sigma$ such that, setting $\sigma=\sigma_\Sigma$, 
\begin{equation}\label{hypstmainth21}
\|\overset{\circ}{A}\|_{L^4(\Sigma)}\leq C_0\sigma^{-1-\delta}, \qquad |\vec z_\Sigma|\leq C_0\sigma^{1-\delta}, \qquad \frac{\sigma}{r_\Sigma}\leq 1+C_0^{-1},
\end{equation}
for some $C_0>0$ large enough. Notice that, for $B_1$, $B_2$ and $B_\textup{cen}$ suitably large $\Sigma$ belongs to $\mathcal B_\sigma(B_1,B_2,B_\textup{cen})$. See also \eqref{remark26}. Suppose that the maximal time of existence of the flow, say $T_\textup{max}$, is finite. Then, by Theorem \ref{proofmainth411} we find that also $\Sigma_{T_\textup{max}}$ belongs to $\mathcal B_\sigma(B_1,B_2,B_\textup{cen})$ and thus, by the regularity theory, we can extend the flow past $T_\textup{max}$, which contradicts the maximality. Thus $T_\textup{max}=\infty$.
\begin{flushright}
$\square$
\end{flushright}
\textbf{Convergence.} From Lemma \ref{gronwallst} we see that $\|H-h\|_{L^2(\Sigma_t)}$ decays exponentially as $t \to +\infty$. Since the derivatives of any order of $H$ are uniformly bounded, interpolation estimates imply that they also decay exponentially. Then Sobolev immersion implies that $\|H-h\|_{L^\infty(\Sigma_t)}$ decays exponentially as well. By the bootstrapping argument described in the paragraph above,  the boundedness of the curvatures and \cite[Lemma 8.2]{huisken1} show that $F(\cdot,t)$ converges to a smooth immersion $F_\infty(\cdot)$. In particular, since $\mathcal H-\hbar\to 0$, the limit surface $\Sigma_\infty:=F_\infty(\Sigma)$ satisfies $\mathcal H\equiv \hbar$. Finally, Theorem \ref{proofmainth411} also shows that the requirements in the definition of $\mathcal{B}_\sigma(B_1,B_2,B_\textup{cen})$ still hold as strict inequalities on $\Sigma_\infty$.
\subsection{CSTMC foliation and centers of mass}
In conclusion, we remark that the computation carried out in the above Lemmas also have some consequences on the center of mass of the foliation we constructed. In this Section, we suppose that the initial datum of our flow is a leaf of Nerz's foliation, as recalled in \eqref{remark26}. In particular, we assume that there exists the CMC-center of mass of Nerz's foliation, i.e. the limit as $s\to\infty$ of the Euclidean barycenters of the foliation $\left\{\Sigma^s\right\}_{s\geq s_0}$ constructed by Nerz  (see \cite{nerz1}). In the following, we will suppose the change of variable $s\longleftrightarrow \sigma$, with $\sigma(s):=\sigma_{\Sigma^s}$. Thus, we have 
\begin{cor} \label{corfin52}
Let $(\M,\overline \g,\overline \K)$ be a \textup{$C_{\frac12+\delta}^2$-asymptotically flat initial data set} which is constrained and with positive ADM-energy $E_\textup{ADM}>0$. Let $\iota^\sigma:\Sigma^\sigma\hookrightarrow \M\setminus \textup{C}$ the inclusion of the family $\left\{\Sigma^\sigma\right\}_{\sigma\geq\sigma_0}$ of \textup{CMC}-surfaces as above and suppose that there exists the \textup{CMC}-center of mass of $\Sigma^\sigma$, i.e. 
\begin{equation}
\vec{\mathcal C}_{\textup{CMC}}:=\lim_{\sigma\to\infty} \fint_{\Sigma^\sigma} \vec x \ d\mu^\sigma,
\end{equation}
where $d\mu^\sigma$ is the 2-dimensional measure induced by $\overline \g$ on $\Sigma^\Sigma$. Consider the \textup{CSTMC} foliation constructed above, and let $\vec z_{\Sigma_\textup{st}^\sigma}$ be the barycenter of $\Sigma_{\textup{st}}^\sigma:=\displaystyle\lim_{t\to\infty} F_t(\Sigma^\sigma)$. 
\begin{enumerate}[label=\textup{(\roman*)}]
\item If $q>\frac2{\frac12+\delta}$ then 
\begin{equation}
\lim_{\sigma\to\infty} \vec z_{\Sigma^\sigma_{\textup{st}}}=\vec{\mathcal C}_{\textup{CMC}}.
\end{equation}
\item If $2\leq q\leq \frac{2}{\frac12+\delta}$, then there exists $C>0$ such that 
\begin{equation}\label{divergentcondition251}
|\vec z_{\Sigma_\textup{st}^\sigma}-\vec z_{\Sigma^\sigma}|\leq C\sigma^{2-\frac{q}2-q\delta}.
\end{equation}
\end{enumerate}
\end{cor}
\begin{proof} Integrating \eqref{barycenterst449} in $[0,t]$ we get 
\begin{equation}\label{eq516}
||\vec z(t)|-|\vec z(0)||\leq C\int_0^t|\Sigma_t|^{-1}\|\mathcal H-\hbar\|_{L^1(\Sigma_t)} \ dt\leq C\sigma^{2-\frac{q}2-q\delta}\left(1-e^{-\frac{E_\textup{ADM}t}{\sigma^3}}\right).
\end{equation}
Since, by construction, $\vec z(0)=\vec z_{\Sigma^\sigma}$ and $\vec z_{\Sigma_\textup{st}^\sigma}:=\lim_{t\to\infty} \vec z(t)$, which exists since the flow converges, letting $t\to \infty$ in \eqref{eq516} we get 
\begin{equation*}
|\vec z_{\Sigma_\textup{st}^\sigma}-\vec z_{\Sigma^\sigma}|\leq C\sigma^{2-\frac{q}2-q\delta}.
\end{equation*}
\end{proof}
\begin{rem}\label{casop2} 
\begin{enumerate}[label=(\roman*)]
\item Since $\frac12+\delta\in \left(\frac12,1\right]$, we have that 
\begin{equation}
2\leq \frac{2}{\frac12+\delta}<4.
\end{equation}
Thus, if $q\geq 4$, the volume preserving spacetime mean curvature flow recovers the center of mass $\vec{\mathcal{C}}_\textup{CMC}$ for every $\delta\in (0,\frac12]$. 
\item For $q=2$, we recover the foliation constructed in \cite{cederbaum}. In this case, the right hand side of equation \eqref{divergentcondition251} is divergent, and the theory developed by Cederbaum and Sakovich in \cite{cederbaum} let us to conclude that $\{\vec z_{\Sigma_\textup{st}^\sigma}\}_{\sigma\geq \sigma_0}$ converges if and only if the correction term converges
\begin{equation}\label{correctterm}
\lim_{r\to\infty}\int_{\mathbb{S}_r^2}\frac{x^i\left(\sum_{k,l}\pi_{kl}x^kx^l\right)^2}{r^3} \ d\mu^e,
\end{equation}
under the additional hypothesis that $|\overline{\K}|\leq \overline c|\vec x|^{-2}$. 
\item Finally, also in the case $q\in \left(2,\frac2{\frac12+\delta}\right]$ equation  \eqref{divergentcondition251} holds with a positive exponent, and thus, in a case in which the CMC-barycenter does not converge, this does not necessarily imply the non convergence of the CSTMC-barycenter. However, differently from point (ii), where the convergence of the limit \eqref{correctterm} allows to deduce a relation between the two barycenters, for a general $q$ we do not know if a similar correction term can be found.
\end{enumerate}
\end{rem}


\begin{thebibliography}{99}
\bibitem{eichmairmots} Andersson, Lars, Michael Eichmair, and Jan Metzger: \emph{Jang’s equation and its applications to marginally trapped surfaces.} Complex Analysis and Dynamical Systems IV: Part 2 (2011): 13-46.
\bibitem{benandrews} Andrews, Ben: \emph{Fully nonlinear parabolic equations in two space variables.} arXiv preprint math/0402235 (2004).
\bibitem{adm} Arnowitt, Richard, and Stanley Deser, and Charles W. Misner: \emph{Coordinate invariance and energy expressions in general relativity.} Physical Review 122.3 (1961): 997-1006.
\bibitem{beigomurchadha} Beig, Robert, and Niall \'O~Murchadha: \emph{The Poincar\'e{} group as the symmetry group of canonical general relativity.} Annals of Physics 147.2 (1987): 463--498.
\bibitem{carfora} Carfora, Mauro, and Annalisa Marzuoli: \emph{Einstein constraints and Ricci flow---a geometrical averaging of initial data sets.} Mathematical Physics Studies, Springer (2023).
\bibitem{bruhat} Choquet-Bruhat, Yvonne: \emph{General relativity and the Einstein equations}, Oxford Mathematical Monographs, Oxford Univ. Press (2009).
\bibitem{cabezasrivas} Cabezas-Rivas, Esther, and Vicente Miquel:  \emph{Volume preserving mean curvature flow in the hyperbolic space}. Indiana University Mathematics Journal 56.5 (2007): 2061-2086.
\bibitem{cederbaum} Cederbaum, Carla, and Anna Sakovich: \emph{On center of mass and foliations by constant spacetime mean curvature surfaces for isolated systems in General Relativity}. Calculus of Variations and Partial Differential Equations 60.6 (2021): 214.
\bibitem{cederbaumnerz} Cederbaum, Carla, and Christopher Nerz: \emph{Explicit Riemannian manifolds with unexpectedly behaving center of mass}. Annales Henri Poincar\'e 16.7 (2015): 1609--1631.
\bibitem{corvinowu} Corvino, Justin, and Haotian Wu: \emph{On the center of mass of isolated systems}. Classical and Quantum Gravity 25.8 (2008): 085008.
\bibitem{delellismuller1} De Lellis, Camillo, and Stefan Müller: \emph{Optimal rigidity estimates for nearly umbilical surfaces.} Journal of Differential Geometry 69.1 (2005): 75-110.
\bibitem{eichmairmots2} Eichmair, Michael, and Lan-Hsuan Huang, and Dan Lee, and Richard Schoen: \emph{The spacetime positive mass theorem in dimensions less than eight.} Journal of the European Mathematical Society 18.1 (2015): 83-121.
\bibitem{eichmairkoerber} Eichmair, Michael, and Thomas Koerber: \emph{Foliations of asymptotically flat 3-manifolds by stable constant mean curvature spheres}. Journal of Diﬀerential Geometry 128.3 (2024): 1037-1083.
\bibitem{gerhardt} Gerhardt, Claus: \emph{Curvature problems}, Series in Geometry and Topology 39, (2006).
\bibitem{glockle} Glöckle, Jonathan: \emph{Spinors and the Dominant Energy Condition for Initial Data Sets}. PhD Thesis. Universität Regensburg, Germany, 2024.
\bibitem{guilisun} Gui, Yaoting, Yuqiao Li, and Jun Sun: \emph{Stability of the area preserving mean curvature flow in asymptotic Schwarzschild space}. Journal of Functional Analysis 289.7 (2025): 111033.
\bibitem{huang} Huang, Lan-Hsuan: \emph{Foliations by stable spheres with constant mean curvature for isolated systems with general asymptotics.} Communications in Mathematical Physics 300.2 (2010): 331-373.
\bibitem{huisken1} Huisken, Gerhard: \emph{Flow by mean curvature of convex surfaces into spheres}. Journal of Differential Geometry 20.1 (1984): 237-266.
\bibitem{huisken86} Huisken, Gerhard: \emph{Contracting convex hypersurfaces in Riemannian manifolds by their mean curvature.} Inventiones Mathematicae 84.3 (1986): 463-480.
\bibitem{huisken1987} Huisken, Gerhard: \emph{The volume preserving mean curvature flow}. Journal f\"ur Reine und Angewandte Mathematik 382
(1987): 35-48.
\bibitem{huiskenpolden} Huisken, Gerhard and Alexander Polden: \emph{Geometric evolution equations for hypersurfaces}, Calculus of variations and geometric evolution problems (Cetraro, 1996), Springer–Verlag, Berlin, 1999, pp. 45–84.
\bibitem{huiskenyau} Huisken, Gerhard, and Shing-Tung Yau: \emph{Definition of center of mass for isolated physical systems and unique foliations by stable spheres with constant mean curvature.} Inventiones Mathematicae 124.1-3 (1996): 281-311.
\bibitem{huiskenwolff} Huisken, Gerhard, and Markus Wolff: \emph{On the evolution of hypersurfaces along their inverse spacetime mean curvature}. arXiv preprint arXiv:2208.05709 (2022).
\bibitem{krwo} Kr\"oncke, Klaus, and Markus Wolff: \emph{Foliations of asymptotically Schwarzschildean lightcones
by surfaces of constant spacetime mean curvature}. arXiv preprint arXiv:2412.17563 (2024).
\bibitem{ladysoloural} Ladyženskaja, Olga A., and Vsevolod Alekseevič Solonnikov: \emph{Linear and quasi-linear equations of parabolic type}. Vol. 23. American Mathematical Soc., 1988.
\bibitem{li} Li, Haozhao: \emph{The volume-preserving mean curvature flow in Euclidean space}. Pacific Journal of Mathematics 243.2 (2009): 331-355.
\bibitem{hawking68} Hawking, Stephen William: 	\emph{Gravitational radiation in an expanding universe.} Journal of Mathematical Physics 9.4 (1968): 598-604.

\bibitem{metzger} Metzger, Jan: \emph{Foliations of asymptotically flat 3-manifolds by 2-surfaces of prescribed mean curvature}. Journal of Differential Geometry 77.2 (2007): 201-236.
\bibitem{miaotam} Miao, Pengzi, and Luen-Fai Tam: \emph{Evaluation of the ADM mass and center of mass via the Ricci tensor}. Proceedings of the American Mathematical Society 144.2 (2016): 753-761.
\bibitem{nerz1} Nerz, Christopher: \emph{Foliations by stable spheres with constant mean curvature for isolated systems without asymptotic symmetry.} Calculus of Variations and Partial Differential Equations 54.2 (2015): 1911-1946.
\bibitem{perez} Perez, Daniel Raoul: \emph{On nearly umbilical hypersurfaces}. PhD Thesis, University of Z\"urich (2011).
\bibitem{vpmcf} Sinestrari, Carlo, and Jacopo Tenan: \emph{Volume preserving mean curvature flow of round surfaces in asymptotically flat spaces.} arXiv preprint arXiv:2501.13091 (2025).
\bibitem{ye} Ye, Rugang: \emph{Foliation by constant mean curvature spheres on asymptotically flat manifolds.} arXiv
preprint dg-ga/9709020 (1997).
\end{thebibliography}
\end{document}